\documentclass[a4paper,10pt]{amsart}
\usepackage{epsf}  
\usepackage{amsfonts, amssymb, amsthm, amsmath}  
\usepackage{hyperref}
\usepackage{tikz}
\usetikzlibrary{matrix,arrows,decorations.pathmorphing}
\usepackage{MnSymbol}
\usepackage{tikz-cd}
\usepackage{enumerate}
\usepackage{graphicx}
\usepackage{psfrag}
\newtheorem{thm}{Theorem}  

\newtheorem{cor}[thm]{Corollary}  
\newtheorem{lemma}[thm]{Lemma}  
\newtheorem{remark}[thm]{Remark}  
\newtheorem{defn}[thm]{Definition}  
\newtheorem{prop}[thm]{Proposition}  
\newtheorem{claim}[thm]{Claim}  
\newtheorem{example}[thm]{Example}  
\newtheorem{alg}[thm]{Algorithm}
\newtheorem{construction}[thm]{Construction}
\numberwithin{thm}{section}  
\def\pf{\noindent\emph{Proof: }}  
\def\stop{\hfill$\square$}

\providecommand{\totl}[1]{\ensuremath{\lceil #1\rceil }}
\providecommand{\totb}[1]{\ensuremath{\underline{ #1}}}

\DeclareMathOperator{\Aut}{Aut}

\newcommand{\ex}{\bold}

\DeclareMathOperator{\expl}{Expl}

\providecommand{\lrb}[1]{\ensuremath{\left(#1\right)}}
\providecommand{\abs}[1]{\left\lvert #1\right\rvert}

\newcommand{\tc}[1]{\check\rvert_{#1}}
\author{Brett Parker   }
\email{brettdparker@gmail.com}  
\thanks{This work was supported by ARC grant DP1093094.}

\title{Tropical enumeration of curves in blowups of $\mathbb CP^{2}$}  
\begin{document}
\maketitle
\begin{abstract} We describe a method for recursively calculating Gromov--Witten invariants of all blowups of the projective plane. This recursive formula is different from the recursive formulas due to  G\"ottsche and Pandharipande in the zero genus case, and Caporaso and Harris in the case of no blowups. We use tropical curves and a recursive computation of Gromov--Witten invariants relative a normal crossing divisor.
\end{abstract}

This paper  computes Gromov--Witten invariants of blowups of the projective plane from recursively computable relative Gromov--Witten invariants. Recursions calculating some such Gromov--Witten  invariants are already known. G\"ottsche and Pandharipande give a recursive formula for zero-genus Gromov--Witten invariants of arbitary blowups of the plane in \cite{GPqc}. Caporaso and Harris  in \cite{CH} show that the  Gromov--Witten invariants of $\mathbb CP^{2}$ relative to a line may be calculated recursively, giving a method for calculating Gromov--Witten invariants of $\mathbb CP^{2}$ of any genus. This is extended by Vakil in  \cite{VakilSurfaces} to a recursive formula for Gromov--Witten invariants of $\mathbb CP^{2}$ blown up at a small number of points, and further extended by Shoval and Shustin, then Brugall\`e  in \cite{SS,Brugalle} to $\mathbb CP^{2}$ blown up at $7$ and $8$ points respectively. Our recursion gives a formula for (arbitary-genus) Gromov--Witten invariants\footnote{We only count rigid curves with no constraints --- the counts of curves constrained to pass through  a collection of points may be recovered by blowing up these points, and counting rigid curves that intersect the new exceptional divisors once.} that is uniform in the number of blowups.

Behind our recursion are exploded manifolds, \cite{iec}, and the tropical gluing formula from \cite{gfgw}. Section \ref{background} contains some necessary background, including a simplified tropical gluing formula, 
\[n_\gamma=\prod_vn_{\gamma_v}\] 
explained and proved in Theorem \ref{simplified gluing}. 

In Section \ref{relative}, we explain our relative Gromov--Witten potential, 
\[ F=\sum_\Gamma \frac{n_\Gamma}{\abs{\Aut\Gamma}}\] 
and how these relative Gromov--Witten invariants $n_\Gamma$ satisfy a recursion induced by relations
\[\overrightarrow ye^F=e^F\overleftarrow y\ .\]
 Section \ref{gf} explains how these relations follow from our simplified gluing formula, and that each side of the above equation counts a  Gromov--Witten invariant using tropical curves with one end constrained; see  Figure \ref{tec10} for a picture of such tropical curves. Our simplified gluing formula  is then applied in Section \ref{reconstruct} to reconstruct our absolute Gromov-Witten invariants from our relative invariants. We count tropical curve such as those in Figure \ref{tec14} to obtain a formula
 \[G=\Psi(F)+\sum_{i} q^{E_{i}}\] 
 first explained in Section \ref{arc}.
  Our reconstruction involves further relative invariants of $\mathbb CP^{2}$ blown up at one point, computed in Section \ref{last} using our simplified gluing formula and  tropical curves  depicted in figures \ref{tec16}, \ref{tec17}, \ref{tec18} and \ref{tec20}.

\

  I have written a Mathematica program that computes these relative Gromov--Witten invariants, and compared the results to known Gromov--Witten invariants of $\mathbb CP^{2}$, and the results of \cite{GPqc,Brugalle,BryanLeung}, however  further work is required  to get these invariants in a closed form. For example, beyond checking up to degree $12$,  I have been unable to reprove the beautiful formulas of Bryan and Leung \cite{BryanLeung} which count curves of genus $g$ passing through $g$ points in  $\mathbb CP^{2}$ blown up at $9$ points.  This formula also can be proved using topological recursion for certain descendant invariants and the symplectic sum formula for Gromov--Witten invariants as in \cite{IP}, however the consequences in our model are less than obvious. The Mathematica program and a link to  a talk with lots of pictures is available on my website:   \href{https://sites.google.com/view/brettparker/home}{https://sites.google.com/view/brettparker/home}. 

\section{Background}\label{background}

\subsection{Stratified $\mathbb Z$--affine spaces}

Our calculation of Gromov--Witten invariants use tropical curves in various  stratified $\mathbb Z$--affine spaces.
\begin{defn} A {\bf $\mathbb Z$--affine map} $\mathbb R^n\longrightarrow \mathbb R^m$ is a map in the form $x\mapsto Ax+c$ where $A$ is an integer matrix and $c\in \mathbb R^m$. A {\bf $\mathbb Z$--affine polytope} $P\subset \mathbb R^n$ is defined by a finite collection of $\mathbb Z$--affine inequalities. In general some of these inequalities may be strict, but say that $P$ is {\bf complete} if it can be defined without  strict inequalities.  A {\bf face} of $P$ is a sub-polytope defined $\alpha=0$ where  $\alpha: P\longrightarrow [0,\infty)$ is $\mathbb Z$--affine, and the {\bf interior} of $P$ is the complement of all proper sub faces. A {\bf complex of $\mathbb Z$--affine polytopes} is a functor from a finite category to the category of $\mathbb Z$--affine polytopes so that each morphism is sent to a $\mathbb Z$--affine inclusion of a face, and each face is the image of a unique morphism. Such a complex has a realisation as a topological space created by gluing the polytopes using the given morphisms; this topological space has $\mathbb Z$--affine strata corresponding to the image of the interior of each polytope --- call this topological space with its $\mathbb Z$--affine strata a {\bf stratified $\mathbb Z$--affine space}.  A {\bf map of stratified $\mathbb Z$--affine spaces} is then a continuous map which sends strata within strata, and which is $\mathbb Z$--affine on strata.
\end{defn}

\begin{figure}[h]
\includegraphics{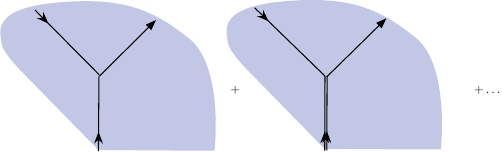}
\caption{}\label{tec13}
\end{figure}

\begin{defn} A {\bf tropical  curve}  is a map of stratified $\mathbb Z$--affine  spaces with a domain consisting of a finite complex of complete $0$--and--$1$--dimensional $\mathbb Z$--affine polytopes. A {\bf vertex} of a tropical curve is a $0$ dimensional stratum, and an ${\bf edge}$ is a $1$--dimensional stratum. An {\bf internal edge} is a bounded $1$--dimensional stratum, and an ${\bf end}$ is an unbounded $1$--dimensional stratum, together with a direction in which that stratum is infinite. 

A tropical curve is {\bf rigid} if it can not be continuously deformed by changing the length of edges and moving vertices within strata. A tropical curve is {\bf rigidly constrained} if a set of its ends is labeled constrained, and it can not be deformed in the above sense without moving at lest one constrained end in a direction independent of the line spanned by that end.

\end{defn}
 For the calculations in this paper, our tropical curves will be in stratified $\mathbb Z$--affine spaces which embed into $\mathbb R^2$, making our tropical curves easy to draw and describe. Examples of tropical curves in different $\mathbb Z$--affine spaces are in Figures \ref{f2}, \ref{tec2}, \ref{tec13}, \ref{tec15}, and \ref{tec16}. 
Our tropical curves will be associated to nodal holomorphic curves, with vertices corresponding to components, ends corresponding to marked points, and internal edges  corresponding to nodes. 

\subsection{Complex manifolds with simple normal crossing divisors}
 A  simple normal crossing divisor $N$ within a complex manifold $M$ is a finite collection of closed, complex codimension 1 submanifolds,   intersecting transversely so that $M$ is covered by holomorphic coordinate charts where $N$ coincides with the coordinate hyperplanes. Such an $(M,N)$ has a natural stratification: a codimension $k$ stratum of $(M,N)$ is a connected component of the set where exactly $k$ of these submanifolds intersect.  Complex manifolds with normal crossing divisors form a category with morphisms   $(M_{1},N_{1})\longrightarrow (M_{2},N_{2})$ consisting of holomorphic maps $M_{1}\longrightarrow M_{2}$ so that each stratum of $(M_{1},N_{1})$ is sent within a  stratum of $(M_{2},N_{2})$,  and the top-dimensional strata of $(M_{1},N_{1})$ are sent within the top-dimensional strata of $(M_{2},N_{2})$. In particular, define a smooth holomorphic curve in $(M,N)$ as a morphism $f:(\Sigma,\{p_i\})\longrightarrow (M,N)$ where $\Sigma$ is a compact and $1$--dimensional. 

We can associate a dual intersection complex to a complex  manifold with a simple\footnote{Complex manifolds with non-simple normal crossing divisors also have dual intersection complexes, but strata of these complexes can have automorphisms, in which case the dual intersection complex is not a stratified $\mathbb Z$--affine space in our sense.} normal crossing divisor: the dual intersection complex is a complex of $\mathbb Z$--affine polytopes which associates a cone, $\sigma_S\equiv [0,\infty)^k$, to each codimension $k$ stratum, and to each inclusion $S\subset \bar S'$ associates a $\mathbb Z$--linear inclusion of $\sigma_{S'}$ as a face of $\sigma_S$.  Maps of complex manifolds with simple normal crossing divisors induce $\mathbb Z$--linear maps of their dual intersection complexes. For example, a map from a  ball around $0$ in  $(\mathbb C,0)^n$ to $(\mathbb C,0)^m$ must be in the form $z\mapsto (z^{\alpha_1}h_1,\dotsc,z^{\alpha_m}h_m)$ where $\alpha_i\in\mathbb N^n$ and $h_i$ are holomorphic and non-vanishing. The induced map on dual intersection complexes $[0,\infty)^n\longrightarrow [0,\infty)^m$ is then $x\mapsto (x\cdot\alpha_1,\dotsc,x\cdot\alpha_m)$. In the case of a smooth curve  $f:(\Sigma,\{p_i\})\longrightarrow (M,N)$, the induced map on dual intersection complexes is a tropical curve. This tropical curve records the contact  of $f$ with the divisor $N$ at each special point $p_i$.

\subsection{Toric manifolds and toric varieties}

A convenient class of examples  are toric manifolds ---  complex manifolds $M$ with a holomorphic $(\mathbb C^*)^n$ action which acts freely on one dense orbit. The toric strata of a toric manifold $M$ are the orbits  of this action. Each codimension $k$ stratum $S$ has a neighborhood isomorphic to $\mathbb C^k\times (\mathbb C^*)^{n-k}$ on which $(\mathbb C^*)^n$ acts, after a change of basis, with the standard action --- this neighborhood consists of all strata with closure containing our codimension $k$ stratum. The dual intersection complex of a toric manifold is called the toric fan, and has extra structure: it naturally embedds into $\mathbb R^n$, which we can identify as a subgroup of $(\mathbb C^*)^n$ using the exponential map.   Each codimension $k$ stratum $S$ determines a $k$--dimensional subset $\sigma^\circ_S\subset\mathbb R^n$ consisting of points $x\in\mathbb R^n$ so that $\lim_{t\to \infty}e^{tx}*p\in S$ for any $p$ in the top dimensional stratum of $M$. The closure, $\sigma_S$, of $\sigma_S^\circ\subset \mathbb R^n$ is a cone, and the set of such   $\sigma_S\subset \mathbb R^n$ is called the fan of $M$. When $S$ is a codimension $k$ stratum of a toric manifold, there is a $\mathbb Z$--linear change of coordinates on $\mathbb R^n$ so that $\sigma_S$ is the positive span of the first $k$ coordinates. In the more general case that $M$ is a toric variety (which for this paper will always mean a normal toric variety), the polytope $ \sigma_S$ is a closed cone not containing any line. In each case, the toric fan is some subdivision of a subset of $\mathbb R^n$, which we consider as the a subgroup of $(\mathbb C^*)^n$ using the exponential map. 

 Toric manifolds and toric varieties are determined by their fans: The standard chart containing a stratum $S$ is $\text{Spec }\mathbb C[  \sigma_S^\vee\cap \mathbb Z^n]$, where $  \sigma_S^\vee\cap \mathbb Z$ is the monoid of $\mathbb Z$--linear functions on $\mathbb R^n$ that are non-negative on $\sigma_S$; the stratum $S$ within this chart is the vanishing locus of the monomials corresponding to $\mathbb Z$--linear functions strictly positive on $\sigma_S^\circ$.   For $S_1\subset   \bar S_2$, $\sigma_{S_2}\subset\sigma _{S_1}$, and the coordinate change map is the natural inclusion $\text{Spec }\mathbb C[  \sigma_{S_2}^\vee\cap \mathbb Z^n]\subset \text{Spec }\mathbb C[  \sigma_{S_1}^\vee\cap \mathbb Z^n]$. More generally,  each equivariant map $f$ of toric varieties (sending the top dimensional stratum within the top dimensional stratum) induces a $\mathbb Z$--linear map of fans, $l$, so that when $f$ sends $S$ within $S'$, $l$ maps $\sigma_S$  within $\sigma_{S'}$. Moreover $f$ in our standard coordinates is some constant times the map $\text{Spec }\mathbb C[  \sigma_{S}^\vee\cap \mathbb Z^n]\longrightarrow \text{Spec }\mathbb C[  \sigma_{S'}^\vee\cap \mathbb Z^{n'}]$ induced by $l$. Conversely, any $\mathbb Z$--linear map of toric fans sending strata within strata induces an equivariant map of toric varieties, determined up to the action of the complex torus on the codomain. 
 
 \subsection{Log geometry }
 
Both complex manifolds with normal crossing divisors and  toric varieties have a natural log structure encoded by the sheaf of holomorphic functions which never vanish within top dimensional strata.  Such functions form a sheaf of monoids which is the log structure sheaf.  Morphisms between these log schemes are the same as the morphisms of toric varieties or complex manifolds with normal crossing divisors described above.

 All log schemes in this paper will be over $\mathbb C$, and log smooth, which means that they are locally isomorphic to an open subset of a  toric variety with its natural log structure (using the analytic  topology). 

\subsection{Exploded manifolds}

\

Holomorphic exploded manifolds generalise complex manifolds with normal crossing divisors, toric varieties, and log smooth log schemes. Exploded manifolds have a  well behaved Gromov--Witten theory, with a degeneration and tropical gluing formula for Gromov--Witten invariants under normal crossing degenerations.  This section lists some neccesary facts about exploded manifolds; for an introduction to exploded manifolds, see one of \cite{iec, scgp, elc}.

\subsubsection{The explosion functor}

\

 The explosion functor, $\expl$, is a functor to the category of holomorphic exploded manifolds from the category of complex manifolds with normal crossing divisors, toric varieties, or log smooth log schemes; for definitions, see  Section 5 of \cite{iec}, Section 2.1 of \cite{scgp},  or Section 4 of \cite{elc}. We define the relative Gromov--Witten invariants of a complex manifold with normal crossing divisors $(M,N)$ to be the Gromov--Witten invariants of $\expl (M,N)$, constructed as in section 8 of \cite{egw}. Alternatively, using log Gromov--Witten theory \cite{acgw,GSlogGW}, or Ionel's method for defining GW invariants relative normal crossing divisors \cite{IonelGW} give the same invariants, as discussed in \cite{elc} and \cite{tropicalIonel}.

\subsubsection{The smooth part and tropical part functors}

\

There are two important functors from the category of exploded manifolds: the smooth part functor, and the tropical part functor. The smooth part, $\totl{\ex B}$, of an (holomorphic) exploded manifold, $\ex B$, is a stratified space with each stratum a smooth (complex) manifold. In particular, $\totl{\expl (M,N)}$ is $(M,N)$,  and the smooth part of the explosion of a log scheme, $M^\dag$,   is the underlying scheme, $M$, with a stratification determined by the log structure. In general,  the tropical part, $\totb{\ex B}$ of an exploded manifold $\ex B$ is a stratified space with each stratum the quotient of some polytope by a group of automorphisms, however this group of automorphisms is  trivial for every exploded manifold appearing in this paper, and $\totb{\ex B}$ is then a stratified $\mathbb Z$--affine space.  For example, $\totb{\expl (M,N)}$ is the dual intersection complex of $(M,N)$, the tropical part of the explosion of a toric variety is its toric fan, and $\totb{\expl M^\dag}$ is the tropicalisation of the log scheme $M^\dag$.

Maps of exploded manifolds induce maps of their smooth and tropical parts, respectively smooth or $\mathbb Z$--affine on each stratum. Given a morphism $f:(M_1,N_1)\longrightarrow (M_2,N_2)$, the smooth part of the explosion of $f$ is just $f$, and the tropical part of the explosion of $f$ is the induced  $\mathbb Z$--linear map on dual intersection complexes. In general, the smooth part of a holomorphic morphism $\tilde f:\expl (M_1,N_1)\longrightarrow \expl (M_2,N_2)$ is a holomorphic map $M_1\longrightarrow M_2$ which respects the stratified structure, but which might not send top dimensional strata within top dimensional strata. The tropical part of $\tilde f$ is $\mathbb Z$--affine on each stratum, and might not be linear. In fact, the tropical part of $\tilde f$ is $\mathbb Z$--linear if and only if the smooth part of $\tilde f$ is a morphisms $\totl{\tilde f}:(M_1,N_1)\longrightarrow (M_2,N_2)$, in which case, $\expl\totl{\tilde f}=\tilde f$. 

\subsubsection{The moduli stack of holomorphic curves}

\

The case of a holomorphic curve in an exploded manifold, $f:\ex C\longrightarrow \ex B$ is of most interest to us; see section 2.3 of \cite{scgp}, section 5 of \cite{elc} or section 8 of \cite{iec} for definitions. The smooth part of $f$ is a nodal holomorphic curve $\totl f$ in the smooth part $\totl{\ex B}$ of $\ex B$. The tropical part of $f$ is a tropical curve $\totb f$ in the tropical part $\totb{\ex B}$ of $\ex B$. The stratifications of $\totl f$ and $\totb f$ are dual; an internal edge of $\totb f$ corresponds to a node of $\totl f$, an end of $\totb f$ corresponds to a marked point of $\totl f$, and a vertex of $\totb f$ corresponds to a component of $\totl{f}$, with nodes and marked points corresponding to the attached internal edges and ends of $\totb{f}$. An end of a holomorphic curve $f$ means a stratum with tropical part an end of $\totl f$, or equivalently, with smooth part a marked point of $\totl f$.

Stable holomorphic curves in a holomorphic exploded manifold $\ex B$  form a stack over the category of exploded manifolds. Moreover, so long as $\ex B$ is complete,\footnote{An exploded manifold $\ex B$ is complete if $\totl{\ex B}$ is compact and $\totb{\ex B}$ is a complex of complete polytopes.}  has a taming form,\footnote{ Taming forms are defined in Section 2 of \cite{cem}. A K\"ahler form on $(M,N)$ gives a taming form on $\expl (M,N)$.} and has tropical part which immerses in $\mathbb R^n$, each connected component of this stack is compact, \cite{cem}, and admits a virtual fundamental class, \cite{vfc}.  As in the case of smooth manifolds, the real dimension of this virtual fundamental class is \begin{equation}\label{vdim} \text{virtual dimension}=2((3-n)(g-1)+ c_1+k)\end{equation} where $n$ is the complex dimension of $\ex B$, and this formula computes the dimension of the component comprised of curves with genus $g$, $k$ ends, and with $c_1$ the first Chern class\footnote{ This Chern class uses the tangent space to $\ex B$, defined in Section 6 of \cite{iec}, or Section 2.2 of \cite{scgp}; the tangent bundle of $\expl (M,N)$ corresponds to the logarithmic tangent space, with holomorphic sections the holomorphic vectorfields tangent to the divisor, so, in the case of the explosion of a smooth curve in $(M,N)$,  $c_1$ is the first Chern class of $M$ evaluated on the curve, minus the degree of intersection of the curve with $N$.} of $\ex B$ evaluated on the curve. Numerical Gromov--Witten invariants can be extracted from the virtual fundamental class by integrating differential forms over it, as in \cite{vfc}. In particular, when the virtual dimension is $0$ we obtain a virtual count of rigid curves by integrating $1$ over the virtual fundamental class.

\subsubsection{Contact data, and counting curves with constrained ends.}

\

The contact data of a smooth holomorphic curve $f$ in $(M,N)$ records the degree of intersection of the curve with components of the divisor $N$ at each special point of the curve, which is also recorded in the tropical part of the curve  $\expl f$ in $\expl (M,N)$. There is an analogous notion of discrete contact data for more general curves in an exploded manifold $\ex B$. 

\begin{defn}An isotopy class of end in $\totb{\ex B}$ is an isotopy class of $\mathbb Z$--affine map $(0,\infty)$. In the exceptional case of an end of a tropical curve with domain $\mathbb R$, the isotopy class of the positive end of a $\mathbb Z$--affine map $\mathbb R\longrightarrow \ex B$ is the isotopy class of the restriction of the map to $(0,\infty)$; otherwise,  the isotopy class of an end of a tropical curve in $\totb{\ex B}$, is the isotopy class of the restriction of the tropical curve to that end.
\end{defn}

In this paper, our exploded manifolds will have tropical part embedded in $\mathbb R^{2}$ so that two ends represent the same isotopy class if and only they have the same derivative, so we can represent isotopy classes of ends in $\totb{\ex B}$ as vectors in some subset of $\mathbb Z^{2}$.

\begin{defn} \label{contact data} {\bf Contact data} $\Gamma$ for curves in an exploded manifold $\ex B$ is a finite set $S_{\Gamma}$, an equivalence relation on $S_{\Gamma}$, and an isotopy class of end in $\totb{\ex B}$ for each element of $S_{\Gamma}$. An {\bf isomorphism of contact data} $\Gamma'\longrightarrow \Gamma$ is a bijective map $S_{\Gamma'}\longrightarrow S_{\Gamma}$ preserving the equivalence relations and the assignment of isotopy classes of ends.  The {\bf contact data of a tropical curve} $\gamma$ in $\totb{\ex B}$ is the set of ends of $\gamma$, along with their isotopy classes, and the equivalence relation joining  two ends if they are in the same connected component of $\gamma$. The contact data of a curve in $\ex B$ is the contact data of its tropical part. A {\bf (tropical) curve with contact data $\Gamma$} is a (tropical) curve and  an isomorphism between $\Gamma$ its contact data, and an isomorphism between (tropical) curves with contact data $\Gamma$ is an isomorphism between the underlying (tropical) curves compatible with the identification of their contact data with $\Gamma$.
\end{defn}

\begin{remark}\label{contact gamma} If $\Gamma$ is contact data for curves in $\expl (M,N)$,  there is a unique rigid tropical curve with contact data $\Gamma$, so we can identify rigid tropical curves in $\totb{\expl (M,N)}$ with contact data. More generally, when $\ex B$ is a cone around a point $0\in \ex B$,  there is a unique tropical curve with contact data $\Gamma$ and with no internal edges, and all vertices sent to $0$. In this case, we identify contact data with such tropical curves.

 When we also have an embedding of $\totb{\ex B}$ into $\mathbb R^{2}$, we will identify contact data for connected curves as a collection of vectors in $\mathbb Z^{2}$, and use the notation $(a,b)\in\Gamma$ to denote such vectors. For curves with multiple connected components, we can encode connected components using an equivalence relation on such a collection of vectors. \end{remark}

There is a moduli stack, $\mathcal M_{\Gamma}$ of curves with contact data $\Gamma$. Because of the extra information encoded in the isomorphism between a curve's contact data and $\Gamma$, this moduli stack is an $\Aut\Gamma$--fold cover of some connected components of the moduli stack of curves in $\ex B$. We can construct an associated virtual fundamental class $[\mathcal M_{\Gamma}]$ as in \cite{vfc}. Usually, we might also want to restrict to components of $\mathcal M_\Gamma$ containing curves with a specified genus and representing some specified homology class, but for the exploded manifolds  in this paper, the homology class represented by a curve is determined by $\Gamma$, and the real virtual dimension of each component is some constant plus $2g$. As explained in \cite{vfc}, we can obtain numerical invariants by pulling back differential forms using evaluation maps, then integrating them over this virtual fundamental class $[\mathcal M_\Gamma]$.

For each end $e\in S_\Gamma$, there are evaluation maps 
\begin{equation}\label{Xe}\begin{tikzcd}\mathcal M_\Gamma \rar{\mathfrak{ev}_e}\ar[bend left]{rr}{ev_e}& \mathcal  X_{e}\rar &\ex X_e\end{tikzcd}\end{equation}
recording the position of the end $e$, and constructed in section 3.4 of \cite{gfgw}. We will not need the precise construction of these evaluation maps, but we will use that, when the derivative of the end is nonzero, $\ex X_e$ is a connected  exploded manifold with complex dimension $\dim \ex X_e=\dim \ex B-1$, and when the derivative of $e$ is $m_e$ times a primitive integral vector,  $\mathcal X_e$ is the stack-quotient of $\ex X_e$ by a trivial $(\mathbb Z/m_e\mathbb Z)$--action. When the derivative of $e$ is zero, $\ex X_e=\ex B$, and $\mathcal X_e$ is the stack-quotient of $\ex B$ by a trivial $\ex T$--action, where $\ex T$ is an exploded manifold with a group structure,  compactifying  $\mathbb C^*$, and described further in section \ref{refinements}.

\begin{defn} Constrained contact data $\Gamma$  is contact data $\Gamma$, with a set $S\subset S_{\Gamma}$ of constrained ends, and for each constrained end $e\in S$, a constraint point $p\in\totb{\ex X_{e}}$. A tropical curve with constrained contact data $\Gamma$ is then a tropical curve $\gamma$ with contact data $\Gamma$ so that for all constrained ends $e\in S$,  $\totb{ev_{e}}(\gamma)=p$.
\end{defn}

\begin{defn} \label{constrained count}
Let $\Gamma$ be constrained contact data. The virtual number of curves with constrained contact data $\Gamma$ is 
\[n_\Gamma:=\int_{[\mathcal M_\Gamma]}\prod_{e\in S}ev_e^*(m_e\theta_e) \]
where the derivative of the end $e$ is $m_e$ times a primitive integral vector, and $\theta_e$ is  Poincar\'e dual to a point in $\ex X_e$  with tropical part the constraint point $p\in \totb{\ex X_{e}}$. 
\end{defn}

The Ponincar\'e dual to a point in $\ex X_e$ is defined in \cite{dre} as the pushforward of $1$ under a map of a point to $\ex X_e$. The integral defining $n_{\Gamma}$ does not depend on the choice of where this point is sent, however choices with different image in $\totb{\ex X_e}$  will lead to different expressions for the same Gromov--Witten invariants, so we call these different constraints, even though they have the same count of constrained curves.

Another way to think of constraining an end $e$ is by picking a map of a point $p\longrightarrow \mathcal X_e$, then defining the constrained moduli stack as $\mathcal M_\Gamma\times_{\mathcal X_e}p$. Using the stack $\mathcal X_e$ instead of the exploded manifold $\ex X_e$ leads to a cleaner gluing formula; see Theorem \ref{simplified gluing}.

%

Let us understand the tropical part of the evaluation map $ev_e$.  If $e$ is nonconstant, $\totb{\ex X_e}$ has strata indexed by the strata of $\totb{\ex B}$ containing some end isotopic to $e$, and is the $\mathbb Z$-affine quotient of these strata with fibers isotopic to $e$.  Moreover,  $\totb{ev_e}$ applied to the tropical part of a curve in $\mathcal M_\Gamma$ is the image of the end labeled by $e$ under this quotient map. So, $\totb{ev_{e}}$ records the position of the end labeled by $e$, where two ends are regarded as having the same position if they eventually coincide.

The Poncare dual to a point in $\ex X_e$ with image $\totb p\in \totb{\ex X_e}$ has support contained in the inverse image of $\totb p$, so the pullback of the Poincar\'e dual to this point to $\mathcal M_\Gamma$ is supported where curves have tropical part with ends constrained to $\totb p$. Constraining an end of a tropical curve in this way constrains the end to eventually lie on a given ray in $\totb{\ex B}$.

\subsubsection{Refinements}\label{refinements}

\
 
A refinement of an exploded manifold $\ex B$ is a complete, bijective submersion $\ex B'\longrightarrow \ex B$.  Each refinement of $\ex B$ has tropical part a subdivison of $\totb{\ex B}$, and any subdivision of $\totb{\ex B}$ into a stratified $\mathbb Z$--affine space uniquely defines a refinement of $\ex B$. Examples include:
\begin{itemize}
\item the explosion of the blowup of $(M,N)$ at a stratum;
\item the explosion of any surjective and generically injective toric map;
\item and the explosion of a logarithmic modification of a log smooth log scheme --- such logarithmic modifications are locally modelled on surjective and generically injective toric maps. 
\end{itemize}
The relationship with subdivisions is familiar in the toric case: Any subdivision of a toric fan into another fan determines  a surjective and generically injective map of toric varieties $M'\longrightarrow M$, and any surjective and generically injective map of toric varieties induces a subdivision of fans.

All compact $n$--dimensional toric varieties have explosions related by refinement. Any two of these exploded manifolds share a common refinement. Moreover, they are all the refinement of a single exploded manifold:  There is a unique $n$--dimensional exploded manifold, $\ex T^n$, with tropical part $\mathbb R^n$. For any compact $n$--dimensional toric manifold $M$, there is a refinement map $\expl M\longrightarrow\ex T^n$ with tropical part the toric fan of $M$ within $\totb{\ex T^n}=\mathbb R^n$. Although, $\ex T$ is not the explosion of a log scheme, it can be regarded as the explosion of the logarithmic torus  $\bf G_\text{log}$; \cite{Glog}.

The moduli stack of holomorphic curves in a refinement,  $\ex B'\longrightarrow \ex B$, is a refinement of the moduli stack of holomorphic curves in $\ex B$. Moreover, Gromov--Witten invariants of $\ex B'$ are equivalent to Gromov--Witten invariants of $\ex B$; a proof of this is  in \cite{egw}, and an analogous result for log Gromov--Witten invariants is proved in \cite{ilgw}. In particular, our counts of curves with constrained ends,  Definition \ref{constrained count}, are unaffected by refining $\ex B$.

Suppose that $\ex B'\longrightarrow \ex B$ is a refinement. Then, $\ex B'\times_{\ex B}\ex B'=\ex B'$. It follows that,  for a cohomology theory of exploded manifolds to have pushforwards that reflect fiber products, the cohomology of $\ex B$ must also contain the cohomology of every refinement $\ex B'$ of $\ex B$. Refined cohomology, defined in \cite{dre}, does just this.  The full tropical gluing formula from \cite{gfgw} uses refined cohomology, but we will only require a simplified version, Theorem \ref{simplified gluing}, which only uses top-dimensional refined cohomology. The top-dimensional refined cohomology of a complete connected exploded manifold is 1 dimensional, as usual, and the refined cohomology class represented by the Poincar\'e dual to a point does not depend on the choice of point. 

\subsection{Decomposition of counts of curves with constrained ends}

\

Counts of rigid curves in an exploded manifold decompose as a sum of well-defined contributions $n_\gamma$ of rigid tropical curves $\gamma$; moreover, our count, $n_\Gamma$, of curves with constrained contact data $\Gamma$ decomposes as a sum over well-defined contributions, $n_\gamma$,  of rigidly constrained tropical curves $\gamma$ with constrained contact data $\Gamma$. In particular, 
\begin{equation}\label{decomposition}n_\Gamma =\sum_\gamma \frac{n_\gamma}{\abs{\Aut \gamma}} \end{equation}
where the sum is over rigidly constrained tropical curves $\gamma$ with constrained contact data $\Gamma$; for interpreting equation (\ref{decomposition}) recall that an automorphism of $\gamma$ is an automorphism of the underlying tropical curve that leaves ends fixed.
  This decomposition follows from Lemma 7.10 of \cite{vfc}, which breaks an integral over the virtual fundamental class up into contributions depending on tropical data, together with Theorem 5.3 and Proposition 5.9 of \cite{evc} which allow us to relate the tropical part of the virtual fundamental class to the moduli space of tropical curves.  An analogous\footnote{When we place our constraint away from the origin in $\totb{\ex X}_e$, some yoga working with log schemes over different bases is required to interpret our constrained curve counts in the log Gromov--Witten setting; nevertheless, with this yoga understood, \cite{acgsdegeneration} implies an analogous decomposition for constrained counts of curves with rational constraints.} decomposition of log Gromov--Witten invariants is proved in \cite{acgsdegeneration}. 
  
\begin{remark}\label{decomposition remark}We will use (\ref{decomposition}) in several contexts. To reconstruct our absolute Gromov--Witten invariants from relative invariants,    we must keep track of the homology class and genus of the curves we are counting. For this, we will need that $n_{\gamma}$ counts the virtual number of curves with an isomorphism of their tropical part with $\gamma$; see \cite{gfgw}, in particular, Lemma 4.8,  the moduli stack $\mathcal M_{[\gamma]}$, evaluation map $ev^{[\gamma]}$, and Gromov--Witten invariant $\eta^{[\gamma]}$. In our context, we will use that $n_{\gamma}$ counts curves with genus and homology class determined by $\gamma$.
\end{remark}

The gluing formula (1) from \cite{gfgw} expresses $n_\gamma$  in terms of local invariants from the vertices of $\gamma$.\footnote{In the setting of log Gromov--Witten invariants,  \cite{rang} proves an analogous gluing formula using the decomposition from \cite{acgsdegeneration} and invariance under logarithmic mondifications from \cite{ilgw}.} This gluing formula uses refined cohomology, but for the exploded manifolds in this paper, we can use a simpler version, Theorem \ref{simplified gluing}.

\subsubsection{Tropical completion}

\

The local invariants from the vertices of a tropical curve $\gamma$ in $\totb {\ex B}$ are constructed using tropical completion. Tropical completion, defined in section 7 of \cite{vfc}, is a canonical way of extending a stratum of $\ex B$ to obtain a complete exploded manifold. For $p\in \totb{\ex B}$, use $\ex B\rvert _{p}$ for the subset of $\ex B$ with tropical part $p$. The tropical completion, $\ex B\tc p$, contains $\ex B\rvert_{p}$ as a dense subset, and the tropical part of $\ex B\tc p$ is a cone around $p$. In particular, we can take the tropical completion of a curve $f:\ex C\longrightarrow \ex B$  at a vertex $v$ of $\totb f$ to obtain its tropical completion 
\[f\tc {v}: \ex C\tc v\longrightarrow \ex B\tc v \]
and this tropically completed curve $f\tc v$ agrees with $f$ restricted to the stratum $\ex C\rvert_{v}$. When the tropical part of $f$ is $\gamma$, the tropical part of $f\tc v$ is a tropical curve $\gamma_{v}$, with a single vertex $v$, and, for each edge leaving $v$, an end with the same derivative.  Our tropical gluing formula for the contribution of $\gamma$ to Gromov--Witten invariants of $\ex B$ uses the Gromov--Witten invariants of $\ex B\tc v$ with the contact data determined by $\gamma_{v}$ as in Remark \ref{contact gamma}. Unless $v$ is sent to a zero--dimensional stratum of $\totb{\ex B}$, $\ex B\tc v$ will not be the explosion of a complex manifold relative a normal crossing divisor, but it will have a refinement that is. 

\subsubsection{Balancing condition of tropical curves}

\

In most tropical geometry, tropical curves satisfy a balancing condition: the sum of the derivatives of edges leaving a vertex is 0. For the tropical part $\totb f$ of a curve in $\ex B$, such a balancing condition will only hold for us at vertices $v$ where $\ex B\tc v$ is a refinement of $\ex T^{2}$.  Then the refinement map  $\ex B\tc v\longrightarrow \ex T^{2}$ gives $\totb{\ex B}\tc v$ as a subdivision of $\mathbb R^{2}$, and we can consider the derivative of the edges of $\totb f$ leaving $v$ as   vectors in $\mathbb Z^{2}$. With this identification, the sum of the derivatives of edges leaving the vertex $v$ is $0$.

At vertices $v$ where $\ex B\tc v$ is not a refinement of $\ex T^{2}$, we will get a weaker balancing condition. In the general case, $\ex B\tc v$ has a refinement equal to some $\expl (M,N)$. Then, the derivatives of edges of $\totb f$ leaving $v$ are the derivatives of the ends of $\totb{f\tc v}$, and this contact data describes the degree of intersection of the smooth curve $\totl {f\tc v}$ in $(M,N)$ with the components of $N$. The balancing condition at $v$ is that this contact data describes the intersection of some homology class with the components of $N$. We will spell out this balancing condition in the specific cases that arise as we progress. It is the usual balancing condition in the case  described above. 

\subsubsection{Tropical gluing formula}

\

The tropical gluing formula, equation (1) from \cite{gfgw}, calculates the contribution of each tropical curve $\gamma$  to Gromov--Witten invariants in terms of a fiber product of  Gromov--Witten invariants of $\ex B\tc v$ with contact data $\gamma_v$ associated to vertices $v$ of $\gamma$.  To describe the gluing formula precisely, all moduli spaces and the spaces over which we take fiber products must be described in terms of exploded manifolds.  For fiber products of exploded manifolds to be reflected correctly using cohomology, it is necessary to use refined cohomology, defined in \cite{dre}, to define the correct Gromov--Witten invariants. To avoid refined cohomology, we shall use a simplified gluing formula that only applies when the following algorithm terminates. 

\begin{alg}\label{alg} Let $\gamma$ be a tropical curve in $\totb{\ex B}$ with some ends constrained, and with each internal edge  attached to two different vertices.
\begin{enumerate}
\item Label any constrained end of $\gamma$  `rigid', and orient it to be incoming.

\item\label{s2} Suppose that a vertex $v$ satisfies the following condition:  the real virtual dimension of the moduli space of curves in $\ex B\tc v$ with unconstrained contact data $\gamma_{v}$ is at least $2(k-1)$ times the number of incoming edges labeled rigid, where $k$ is the complex dimension of $\ex B$.  Then label the remaining edges leaving $v$ as rigid, and orient them away from $v$.

\item Repeat step \ref{s2} until all vertices satisfy the condition from Step \ref{s2}.

\end{enumerate}
\end{alg}

\begin{thm}[Simplified gluing formula]\label{simplified gluing}Let $\gamma$ be a tropical curve in $\totb{\ex B}$ with constrained contact data $\Gamma$.
If  Algorithm \ref{alg} terminates with all edges of $\gamma$ labeled rigid, 
\[  n_{\gamma}=\prod_{v}n_{\gamma_{v}} \]
where $n_{\gamma}$ appears in equation (\ref{decomposition}), and $ n_{\gamma_{v}}$ is the virtual number of curves in $\ex B\tc v$ with  constrained contact data specified by $\gamma_{v}$ with incoming edges constrained and outgoing edges unconstrained. 

\end{thm}

\pf The gluing formula,   Equation (1) from \cite{gfgw}, is a combination of Theorem 4.7,
\begin{equation}\label{gf1}\eta^{[\gamma]}=k_{\gamma}\Delta^{*}\prod_{v}\eta^{[\gamma_{v}]}\ ,\end{equation}
and Lemma 4.8,
\begin{equation}\label{gf2}\eta\tc{\gamma}=\frac1 {\abs{\Aut \gamma}}i^{[\gamma]}_{!}\eta^{[\gamma]} \ .
\end{equation}
%
The form $\eta \tc\gamma$ is\footnote{In \cite{gfgw}, $\eta$ includes dummy variables $q$ and $\hbar$ whose exponents record homology and Euler class respectively. It turns out that such information is encoded in the the tropical part of curves in this paper, so, for this paper, we set $q$ and $\hbar$ to $1$ when we define $\eta\tc \gamma$.} the pushforward of the tropical completion of the virtual fundamental class at $\gamma$, where the pushforward uses the tropical completion of the  evaluation map $\prod_{e\in S_\Gamma}ev_e$. In particular, using notation from Definition \ref{constrained count}, the contribution of $\gamma$ to $n_\Gamma$ is 
%
\[\int_{\prod_{e\in S_{\Gamma}}\ex X_{e}\tc{\totb{ev_{e}}\gamma}}\eta\tc\gamma\wedge \prod_{e\in S}m_{e}\theta_{e}\tc{\totb{ev_{e}}\gamma}\ .\] 

The form $\theta_{e}\tc{\totb{ev_{e}}\gamma}$ vanishes unless $\gamma$ has the end $e$ constrained as specified by the constrained contact data $\Gamma$, in which case, $\theta_{e}\tc{\totb{ev_{e}}\gamma}$ is the Poincare dual to a point in $\ex X_{e}\tc{\totb{ev_{e}}\gamma}$.  Moreover, the integral above vanishes if  $\gamma$ is not rigidly constrained. So only rigidly constrained tropical curves with constrained contact data $\Gamma$ contribute to $n_{\Gamma}$. For the remainder of this proof, we assume that $\gamma$ has constrained contact data $\Gamma$.

Now,  associate an exploded manifold, $\ex X_e$, to each edge $e$ of $\gamma$:  for each vertex $v$ connected to $e$, $\ex X_{e}$ is the codomain of the evaluation map $ev_{e}$ from (\ref{Xe}) at the corresponding end of $\gamma_{v}$. For $e$ an end of $\gamma$, this $\ex X_{e}$ is a tropical completion, at $\totb{ev_{e}}\gamma$, of the $\ex X_{e}$ above associated to $e\in S_{\Gamma}$. 

Use $\text{ied}\gamma$ for the internal edges of $\gamma$, $\text{End}\gamma$ for the set of ends of $\gamma$, and $S\subset \text{End}\gamma$ for the set of constrained ends. The maps $i^{[\gamma]}$ and $\Delta$ in (\ref{gf1}) and (\ref{gf2}) are the obvious projection and diagonal inclusion below.  
\[\begin{tikzcd} \prod_{e\in \text{End}\gamma }\ex X_e&\lar{i^{[\gamma]}} \prod_{e\in \text{End}\gamma} \ex X_e\times \prod_{e\in \text{ied}\gamma} \ex X_e\rar{\Delta} & \prod_{e\in \text{End}\gamma} \ex X_e\times \prod_{e\in \text{ied}\gamma} \ex X^2_e\dar{=}
\\ & & \uar\prod_v\prod_{e\in \text{End}\gamma_v}\ex X_e \end{tikzcd} \]
We can define $n_\gamma$ using the form $\eta^{[\gamma]}$, defined in \cite{gfgw} before Theorem 4.7, but as with $\eta\tc \gamma$, we set $q$ and $\hbar$ to $1$ because $\gamma$ already keeps track of homology class and Euler characteristic for us. 
\begin{equation}\label{ngammadef}n_\gamma:=\int_{\prod_{e\in \text{End}\gamma\cup\text{ied}\gamma} \ex X_e}\eta^{[\gamma]}\prod_{e\in S}m_{e}\theta_{e}\end{equation}
where now, for each edge $e$ of $\gamma$,  the form $\theta_{e}$ is the Poincare dual to a point in $\ex X_{e}$, and $m_{e}$ is the multiplicity of $e$.
Using (\ref{gf2}), the contribution of $\gamma$ to $n_{\Gamma}$ is 
\[\frac 1{\abs{\Aut \gamma}}\int_{\prod_{e\in \text{End}\gamma\cup\text{ied}\gamma} \ex X_e}\eta^{[\gamma]}\prod_{e\in S}m_{e}\theta_{e}=\frac 1{\abs{\Aut\gamma}}n_{\gamma}\]
\

Now assume that Algorithm \ref{alg} terminates with all vertices satisfying the condition from Step \ref{s2}, and let $I_v$ be the set of incoming edges of $\gamma_v$.  The form $\eta^{[\gamma_{v}]}$ is the pushforward of the virtual fundamental class of curves with contact data $\gamma_{v}$ in $\ex B\tc v$. (It has a similar definition to $\eta^{[\gamma]}$, but a tropical completion complicating the definition of $\eta^{[\gamma]}$ does nothing and is not required for defining $\eta^{[\gamma_v]}$.) The count of curves with contact data $\gamma_{v}$, and ends in $I_{v}$ constrained, is then
\[n_{\gamma_v}=\int_{\prod_{e\in \text{End}\gamma_v}\text X_e}\eta^{[\gamma_v]}\wedge\prod_{e\in I_v }m_e\theta_e\ .\]
Using the condition from Step \ref{s2} of Algorithm \ref{alg}, we can reformulate the above as
\begin{equation}\label{spf}(\pi_{v})_{!}\lrb{\eta^{[\gamma_v]}\wedge\prod_{e\in I_v }m_e\theta_e}=n_{\gamma_{v}}\prod_{e\in \text{End}\gamma_{v}\setminus I_{v}}\theta_{e} +\text{lower dimensional terms,}\end{equation}
where $\pi_{v}$ is the obvious projection
\[\pi_{v}:\prod_{e\in \text{End}\gamma_{v}}\ex X_{e}\longrightarrow \prod_{e\in \text{End}\gamma_{v}\setminus I_{v}}\ex X_{e}\ .\]

The constant $k_\gamma$ from (\ref{gf1}) is  $k_{\gamma}:=\prod_{e\in\text{ied}\gamma}m_{e}$. Equation (\ref{gf1}) then implies
\begin{equation}\label{ngamma}n_{\gamma}=\int_{\prod_{e\in \text{End}\gamma\cup\text{ied}\gamma}\ex X_e}\lrb{\prod_{e\in \text{ied}\gamma}m_{e}}\prod_v\eta^{[\gamma_v]}\wedge  \prod_{e\in S}m_{e}\theta_e\ .\end{equation}
We break the integral (\ref{ngamma}) into a series of pushforwards. Let $O:=\text{End}\gamma \setminus S$, and note that each edge of $\gamma$ is either in $O$, or is in $I_{v}$ for a unique vertex $v$ of $\gamma$.  Let $v_{1},\dotsc, v_{n}$ be the vertices of $\gamma$, ordered so that Step \ref{s2} of Algorithm \ref{alg}  applies to the vertex $v_{i}$ after it has been applied to $v_{i+1}$. Let $\pi_{v_{i}}$ be the obvious projection collapsing $X_{e}$ for $e\in I_{v_{i}}$:
\[\pi_{v_{i}}:\prod_{e\in O\bigcup_{j\leq i}I_{v_{j}}}\ex X_{e}\longrightarrow \prod_{e\in O\bigcup_{j< i}I_{v_{j}}}\ex X_{e}\]
\begin{claim}\label{pf}
\[(\pi_{v_{i}})_{!}\lrb{\prod_{e\in I_{v_{i}}}m_{e}\prod_{j\geq i}\eta^{[\gamma_{j}]}\prod_{e\in S_{i}}\theta_{e}}=n_{\gamma_{v_{i}}}\prod_{j> i}\eta^{[\gamma_{j}]}\prod_{e\in S_{i-1}}\theta_{e}+\text{lower dimensional terms}\]
\[\text{where }S_{n}:=S \text{ and } S_{i-1}:=S_{i}\cup \text{End}\gamma_{v_{i}}\setminus I_{v_{i}}\ .\]

\end{claim}

To prove Claim \ref{pf}, note that we are integrating over fibers $\prod_{e\in I_{v_{i}}}X_{e}$. These incoming edges to $v_{i}$ are not attached to $v_{j}$ for any $j>v_{i}$, so $\eta^{[\gamma_{v_{j}}]}$ vanishes on these fibers. Moreover, $\theta_{e}$ vanishes on these fibers unless $e\in I_{v_{i}}$. Therefore, Claim \ref{pf} follows from equation (\ref{spf}).

Applying Claim \ref{pf} repeatedly to the integrand of (\ref{ngamma}), and noting that the condition from Algorithm \ref{alg} ensures that the lower dimensional terms from Claim \ref{pf} still push forward to have lower dimension than the terms we are interested in, we get
\[\begin{split}(\pi_{1})_{!}\circ \dotsb\circ(\pi_{n})_{!}\lrb{\lrb{\prod_{e\in S\cup\text{ied}\gamma}m_{e}}\prod_v\eta^{[\gamma_v]}\wedge  \prod_{e\in S}\theta_e}=&\prod_{i=1}^{n}n_{\gamma_{v_{i}}}\prod_{e\in O}\theta_{e}\\ &+\text{lower dimensional terms.}\end{split}\]
Then, integrating  over $\prod_{e\in O}\ex X_{e}$ gives our required formula.
\[n_{\gamma}=\prod_{i=1}^{n}n_{\gamma_{v_{i}}}\]

\stop

We also will need the following easy consequence of the full gluing formula. Let $\gamma$ be a rigidly constrained curve in $\totb{\ex B}$, where $\ex B$ is a $2$--complex-dimensional exploded manifold. Define $d_v$  as the minimal real dimension of a nonempty component of the virtual fundamental class of the moduli stack of curves with contact data $\gamma_v$, or set $d_v=\infty$ if this virtual fundamental class is empty. Let the number of constrained ends of $\gamma$ be $\abs S$, and let the number of internal edges of $\gamma$ be $\abs{\text{ied}\gamma}$, and the valence of $v$ be $\abs{\text{End}\gamma_v}$. 
\begin{lemma}\label{dimension vanishing} The contribution $n_\gamma$ of a rigidly constrained tropical curve to Gromov--Witten invariants of a $2$--complex-dimensional exploded manifold $\ex B$ vanishes if either of the following conditions hold:
\item for any vertex $v$ of $\gamma$
\[d_v>\abs{\text{End}\gamma_v} \ ;\]
\item or \[\sum_vd_v -\abs{S}-\abs{\text{ied}\gamma}>0\]
\end{lemma}

\pf In our gluing formula, (\ref{gf1}), the exploded manifolds $\ex X_e$ have real dimension $2$, and the form $\eta^{[\gamma_v]}$ has dimension at most $2(\abs{\text{End}\gamma_v}-d_v)$, so vanishes if this is negative. In this case, (\ref{gf1})  implies that $n_\gamma=0$. Moreover, (\ref{gf1}) implies that $\eta^{[\gamma]}$ has codimension at least $2(\sum_vd_v-\abs{\text{ied}\gamma})$, so the integral (\ref{ngammadef}) defining $n_\gamma$ vanishes if $\sum_vd_v -\abs{\text{ied}\gamma}>\abs S$.

\stop

\section{The relative Gromov--Witten invariants}\label{relative}

%

\begin{figure}[h]\label{f1}
\includegraphics{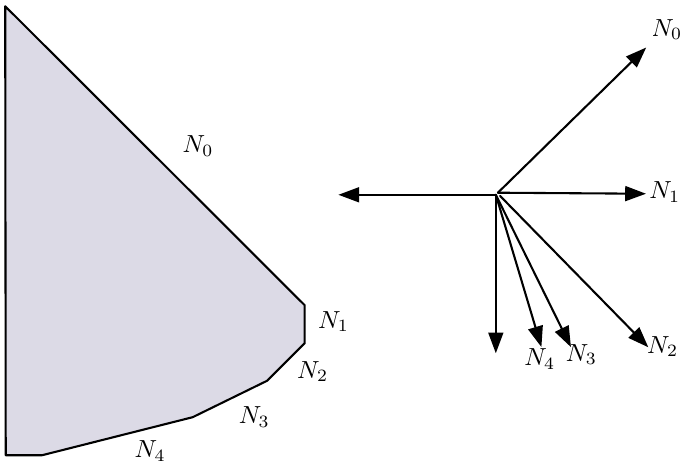}\caption{The moment polytope, and dual fan of $M$}
\end{figure}

Our formula for the Gromov--Witten invariants of blowups $\mathbb CP^2$ involves relative invariants of  $(M,N)$, where $M$ is a toric blowup of $\mathbb CP^2$, and $N$ is a union of toric boundary strata of $M$.
In particular, $M$ has toric fan the span of  $(-1,0)$, $(0,-1)$, $(1,1)$, and $(1,1-k)$ for $k=1,\dotsc n$, and $N=\bigcup_{k=0}^n N_k$ where $N_k$ is the toric stratum of $M$ corresponding to $(1,1-k)$ for $k=0,\dotsc, n$. Figure \ref{f1} shows  a symplectic moment polytope and dual fan  for one such $M$. 

Our relative Gromov--Witten invariants concern holomorphic curves in $M$ with specified contact data with the divisor $N$. We encode such contact data using integral vectors in the non-negative  span of $(1,1)$ and $(1,1-n)$ as follows: A vector 
$(d,d-id)$ indicates a point of contact order $d$ with $N_{i}$, and a vector,  $a(1,1-i)+b(1,1-i-1)$ where $a$ and $b$ are positive integers, indicates  a point  sent to $N_{i}\cap N_{i+1}$, where the curve is required to have contact order $a$ to $N_{i}$ and $b$ to $N_{i+1}$. Note that this contact data determines the homology class represented by a curve.

\begin{figure}[h]
\includegraphics{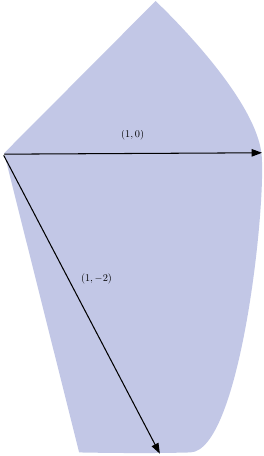}
\caption{ A rigid tropical curve in the tropical part of $\expl(M,N)$, or $\totb {\ex M}$.}\label{f2}
\end{figure}

We will represent the dual intersection complex of $(M,N)$, or the tropical part of $\expl (M,N)$ as the non-negative span of $(1,1),\dotsc,(1,1-n)$.  If $f$ has contact data encoded by a set, $\Gamma$, of integral vectors, the tropical part of $\expl f$ is a rigid tropical curve with a single vertex, sent to $0$, and an end leaving $0$ with derivative $(a,b)$ for each vector, $(a,b)$ in $\Gamma$. One such rigid tropical curve is depicted in Figure \ref{f2}.  More generally, the tropical part of any stable holomorphic curve in $\expl(M,N)$ is a tropical curve with a finite set of ends; the derivatives of the tropical curve at these ends is a set, $\Gamma$, of integral vectors in the non-negative span of $(1,1),\dotsc, (1,1-n)$, and we call $\Gamma$ the contact data of the holomorphic curve. 

%

We will count rigid holomorphic curves with unconstrained contact data $\Gamma$.
%

\begin{defn} Let $n_\Gamma$ be the virtual number of rigid holomorphic curves in $\expl(M,N)$ with contact data $\Gamma$, so
\[n_\Gamma=\int_{[\mathcal M_{\Gamma}]}1\ .\]


\end{defn}
%


\begin{defn}\label{Fndef} The relative Gromov--Witten potential of $(M,N)$ is
\[F_{n}:=\sum_{\Gamma}\frac {n_{\Gamma}}{\abs{\Aut \Gamma}}\Gamma \ .\]
\end{defn}
This relative Gromov--Witten potential takes values in the vectorspace with basis labeled by unordered finite sets $\Gamma$ of integral vectors in the nonnegative span of $(1,1)$ and $(1,1-n)$. We identify such $\Gamma$ with the set of connected rigid tropical curves in the tropical part of $\expl{(M,N)}$, or the contact data for connected curves in $\expl(M,N)$ as in Definition \ref{contact data}.

\subsection{Absolute Gromov--Witten invariants from relative invariants.} \label{arc}
	
	\
	
The purpose of this paper is to compute the absolute Gromov--Witten invariants of blowups of $\mathbb CP^2$. We use the  construction of section 8 in \cite{egw} to define these Gromov--Witten invariants so that we can  apply a tropical gluing formula from \cite{gfgw}.
	
\begin{defn} \label{Gndef} Let $n_{g,\beta}$ be the virtual number of rigid curves in the blowup of $\mathbb CP^2$ at $n$ points, with genus $g$ and representing the homology class $\beta$. The absolute Gromov--Witten  potential is 
\[G_{n}:=\sum n_{g,\beta}x^{g-1}q^{\beta} \ .\]
\end{defn}
	
 We now describe how to obtain the absolute Gromov--Witten potential $G_n$ from the relative potential, $F_n$. Our result, proved in Proposition \ref{GnFn proof}, takes the form
\begin{equation}\label{GnFn}G_{n}=\Psi(F_{n})+\sum_{i=1}^{n}x^{-1}q^{E_{i}} \ ,\end{equation}
with terms explained below.

	The second homology of our blowup of $\mathbb CP^2$ is generated by $H$ and $E_1,\dotsc,E_n$, where  $E_i$ is the homology class of the $i$th exceptional sphere,  and $H$ is  the homology class represented by the pullback of a line. The terms $x^{-1} q^{E_i}$ above represent these $n$ exceptional spheres. 
 
 To define $\Psi(F_n)$, we first define $\Psi(\Gamma)$. For  $\Gamma=\{(1,1-m_{1}),\dotsc,(1,1-m_{k})\}$, define
\begin{equation}\label{Gamma contribution}\Psi(\Gamma):= \prod_{i=1}^{k}\lrb{x^{m_{i}-3}q^{H}\sum_{j_{1}<\dotsb<j_{m_{i}}}q^{-\sum_{l=1}^{m_{i}}E_{j_{l}}}}\end{equation}
and define $\Psi(\Gamma)=0$ if $\Gamma$ is not in the above form. The multiplication above is the $\mathbb R$--linear extension of the rule $(x^{a_1}q^{\beta_1})(x^{a_2}q^{\beta_2}):=x^{a_1+a_2}q^{\beta_1+\beta_2}$.   Note that the sum on the right is the $m_{i}$-th elementary symmetric polynomial $\sigma_{m_{i}}$ in the variables $\{q^{-E_{1}},\dotsc,q^{-E_{n}}\}$, so the above formula is succinctly written as
\[\Psi(\{(1,1-m_{1}),\dotsc,(1,1-m_{k})\})=\prod_{i=1}^{k}x^{m_{i}-3}q^{H}\sigma_{m_{i}}\ .\]

We complete our description of the terms in equation (\ref{GnFn}) by extending the definition of $\Psi$ to be $\mathbb R$--linear, so 
\[\Psi(F_n):=\sum_\Gamma \frac{n_\Gamma}{\abs{\Aut\Gamma}} \Psi(\Gamma)\ .\]
 Although this is an infinite sum, there are only a finite number of $\Gamma$ contributing to the coefficient of $x^{g-1}q^{\beta}$, so it is well defined.   

By allowing the number of blowups to approach infinity, we may write a formula that is uniform in the number of blowups. Let $A$ be the set of vectors in the form $(a,b)\subset \mathbb Z^2$ with $a>0$ and $b\leq a$, and let $S(A)$ be the set of finite, unordered collections of vectors in $A$. Let 
\[F=\sum_{\Gamma\in S(A)}\frac{n_{\Gamma}}{\abs{\Aut\Gamma}}\Gamma\]
where $n_{\Gamma}$ is the corresponding relative Gromov--Witten invariant of $(M,N)$ where the number $n$ of blowups used\footnote{So long as $n$ is large enough for $n_{\Gamma}$ to be defined, the value of $n_{\Gamma}$ does not depend on $n$.} is enough that the vectors in $\Gamma$ are in the non-negative span of $(1,1)$ and $(1,1-n)$.  Then
\begin{equation}\label{translation}G=\Psi(F)+\sum_{i} q^{E_{i}}\end{equation}
is a generating function representing the rigid curves in all possible blowups of $\mathbb CP^{2}$.

\subsection {Recursive calculation of relative Gromov--Witten potential.}\label{recursive calculation}

\

We shall now describe how to recursively compute our relative Gromov--Witten potential $F$. Our recursion is simplified by using $e^{F}$, the generating function counting possibly disconnected curves. Use the convention that $(n_{1}\Gamma_{1})(n_{2}\Gamma_{2}):=(n_{1}n_{2})\Gamma_{1}\coprod\Gamma_{2}$. Then  write the generating function for possibly disconnected curves as
\[e^{F}=\sum_{\Gamma\in SS(A)} \frac{n_{\Gamma}}{\abs{\Aut\Gamma}}\Gamma\]
where the sum is now over the set, $SS(A)$, of unordered, finite\footnote{We include the case $\Gamma=\emptyset$. There is a unique empty curve, so  $n_{\emptyset}=1$.} collections $\coprod_i\Gamma_i$ of  $\Gamma_i \in S(A)$. Such a $\coprod_i\Gamma_i\in SS(A)$ corresponds to a rigid tropical curve, with connected components the rigid tropical curves corresponding to $\Gamma_i$, and gives contact data, as in Remark \ref{contact gamma}, for curves in $\expl(M,N)$ when enough blowups are used so that $\Gamma$ is contact data for a tropical curve in $\totb{\expl(M,N)}$.  The count of disconnected curves is easily read off from counts of connected curves using  the formula  $n_{\coprod\Gamma_{i}}=\prod_{i}n_{\Gamma_{i}}$.  An alternate description of $\Gamma\in SS(A)$ is as a finite set of integral vectors $(a,b)$ so that $a>0$ and $b\leq a$ along with an equivalence relation --- where we say some of these vectors are connected to each other. For our recursive calculation, we will need the following two gradings on the set $SS(A)$ of such $\Gamma$: 
\begin{equation}\label{degref}\text{Degree: } \deg\Gamma:=\sum_{(a,b)\in\Gamma}a\end{equation}
\[\text{Euler Characteristic: } \chi(\Gamma):= \sum_{(a,b)\in \Gamma}(1+2a+2b)  \]
\begin{remark}\label{chi bound} If a connected curve has genus $g$ and $k$ ends, we say its Euler characteristic is $2-2g-k$, and more generally the Euler characteristic of a curve is the sum of its components Euler characteristics. Lemma \ref{deg and vdim} implies that a possibly disconnected  holomorphic curve in $(M,N)$ with contact data $\Gamma$ will have a $(\deg\Gamma)$--fold  intersection with a line from $\mathbb CP^2$, and if it is rigid, it will have Euler characteristic $\chi(\Gamma)$. In particular, because each connected component of a curve has Euler characteristic bounded by $2$, Lemma \ref{deg and vdim} implies that $n_\Gamma$ is zero if $\chi(\Gamma)$ is greater than twice the number of connected components of $\Gamma$. 
\end{remark}

Our recursion will use that the set of $\Gamma\in SS(A)$ with degree bounded above and Euler characteristic bounded below is finite.

The generating function $e^{F}$ is determined by the relations
\[e^{F}\overleftarrow {y}=\overrightarrow{y}e^{F}\ ,\]
where  $y=(-a,-b)$ with $a>0$ and $b<a$,  and  $\overleftarrow y$ and $\overrightarrow {y}$  are linear maps from the vectorspace with basis $SS(A)$ to the vectorspace with basis labeled by $SS(A)_{y}\subset SS(A\cup\{y\})$, the subset comprised of $\Gamma\in  SS(A\cup\{y\})$ containing a unique vector equal to $y$. We also interpret $\Gamma\in SS(A)_y$ as a rigid, possibly disconnected tropical curve with one special end of derivative $-y$ labeled constrained, and all other ends labeled unconstrained. Then we can use Definition \ref{constrained count} to define $n_\Gamma$ as the count of holomorphic curves with contact data $\Gamma$ and this special end constrained. In fact, we will show 
%
%
\begin{equation}\label{free edge}e^{F}\overleftarrow y=\sum_{\Gamma\in \SS(A)_y}\frac{n_{\Gamma}}{\abs{\Aut\Gamma}}\Gamma=\overrightarrow y e^{F}\ \end{equation}
 and show  that $e^{F}\overleftarrow y$ and $\overrightarrow y e^{F}$ calculate $n_{\Gamma}$ by constraining the special end in  two different positions. 

Let us make a combinatorial definition for the linear maps $\overleftarrow y$ and $\overrightarrow y$. Use the notation
\[(a,b)\wedge(c,d):=ad-bc \ .\]
Let $X=\{v_1,\dotsc,v_n\}$ be a set of vectors in $\mathbb Z^2$, ordered so that $v_i\wedge v_j\geq 0$ if $i<j$.
Define
\[n_{X\overleftarrow y}:=\prod_{i=1}^n\max\left\{0, v_i\wedge\lrb{ y+\sum_{j=i+1}^nv_j}\right\}\text{ and }n_{\overrightarrow y X}:=\prod_{i=1}^n\max\left\{0, \lrb{ y+\sum_{j=1}^{i-1}v_j}\wedge v_i\right\}\ .\]
Suppose now that $X$ is a subset of the set of vectors in $\Gamma:=\{\Gamma_1,\dotsc,\Gamma_k\}\subset SS(A)$. Define $\Gamma_{X,y}$ and $\Gamma'_{X,y}$ in $ SS(\mathbb Z^2)$ as follows:
\begin{equation}\label{gammaX}\begin{split}\Gamma_{X,y}&:= \lrb{\{y,y+\sum_{v\in X} v\}\cup \bigcup_{\Gamma_i\cap X\neq \emptyset}\Gamma_i\setminus X}\coprod_{\Gamma_j\cap X=\emptyset}\Gamma_j  
\\ \Gamma'_{X,y}&:= \lrb{\{y\}\cup \bigcup_{\Gamma_i\cap X\neq \emptyset}\Gamma_i\setminus X}\coprod_{\Gamma_j\cap X=\emptyset}\Gamma_j\end{split}  \end{equation}
In words, $\Gamma_{X,y}$ is obtained by combining all sets of vectors $\Gamma_{i}$ containing vectors in $X$, and replacing  the vectors in $X$ with $\{y, y+\sum_{v\in X}v\}$, whereas $\Gamma'_{X,y}$  combines all sets of vectors containing vectors in $X$ then replacing the vectors in $X$ with $y$. We will need $\Gamma_{X,y}$ in the case that it is in $SS(A)_y$, so $y+\sum_{v\in X}v=(a,b)$ with $a>0$ and $b\leq a$, and will use $\Gamma'_{X,y}$ in the case that $y+\sum_{v\in X}v$ is $(-1,0)$ or $(0,-1)$. In both cases, $\Gamma_{X,y}$ and $\Gamma'_{X,y}$ are in $SS(A)_y$.

 In section \ref{tropical description}, we describe rigidly constrained tropical curves with contact data $\Gamma_{X,y}$ and $\Gamma'_{X,y}$. These tropical curves can be described using an incoming constrained edge with derivative $y$ interacting with the vectors in $X$. We shall see that these tropical curves contribute $n_{X\overleftarrow y}$ or $n_{\overrightarrow yX}$ to Gromov--Witten invariants.

Use the notation $X\subset \Gamma$ for a set of vectors in $\Gamma$. Then define
\begin{equation}\label{overleftarrowdef}\Gamma\overleftarrow  y:=\sum_{X\subset \Gamma\vert \Gamma_{X,y}\in SS(A)_y}n_{X\overleftarrow y}\Gamma_{X,y}+\sum_{X\subset \Gamma\vert y+\sum_{v\in X}v=(-1,0)}n_{X\overleftarrow y}\Gamma'_{X,y}\end{equation}
and
\begin{equation}\overrightarrow y\Gamma:=\sum_{X\subset \Gamma\vert \Gamma_{X,y}\in SS(A)_y}n_{\overrightarrow yX}\Gamma_{X,y}+\sum_{X\subset \Gamma\vert y+\sum_{v\in X}v=(0,-1)}n_{\overrightarrow yX}\Gamma'_{X,y}\end{equation}
and finally, 
\begin{equation}\label{left action}\overrightarrow y e^F:=\sum_{\Gamma\in SS(A)}\frac{n_\Gamma}{\abs{\Aut\Gamma}}\overrightarrow y \Gamma\end{equation}
and
\begin{equation}\label{right action} e^F\overleftarrow y:=\sum_{\Gamma\in SS(A)}\frac{n_\Gamma}{\abs{\Aut\Gamma}} \Gamma\overleftarrow y\end{equation}

Below, we will check that the infinite sums in equations (\ref{left action}) and (\ref{right action}) are well defined. For $\Gamma\in SS(A)_{(-y_1,-y_2)}$, define
\[\deg\Gamma:=2y_1+\sum_{(a,b)\in \Gamma}a\]
and
\[\chi(\Gamma):=4y_1+4y_2-2+\sum_{(a,b)\in \Gamma}(1+2a+2b) \ .\]
Then, recalling equations (\ref{degref}) and (\ref{gammaX}) gives the following identities. 
\begin{equation}\label{degchi}\begin{split}\deg \Gamma_{X,y} &=\deg\Gamma 
\\  \chi(\Gamma_{X,y}) &=\chi(\Gamma)-\abs{X}
\\ \deg\Gamma'_{X,y}&=\begin{cases} \deg\Gamma \ \ \ \ \ \ \ \text{ if $y+\sum_{v\in X}v=(0,-1)$} \\ \deg\Gamma+1\ \ \ \text{ if $y+\sum_{v\in X}v=(-1,0)$}\end{cases}
\\ \chi(\Gamma'_{X,y})&=\chi(\Gamma)+1-\abs X\end{split}\end{equation}
Because the number of vectors in $\Gamma\in SS(A)$ is bounded by $\deg\Gamma$, the number, $\abs X$, of vectors in $X\subset \Gamma$  is bounded by $\deg\Gamma$, so we get the following bounds on $\chi(\Gamma)$ in terms of $\Gamma_{X,y}$ and $\Gamma'_{X,y}$. 
\[\chi(\Gamma_{X,y})\leq \chi(\Gamma)\leq \chi(\Gamma_{X,y})+\deg\Gamma_{X,y}\]
\[\chi(\Gamma'_{X,y})-1\leq \chi(\Gamma)\leq \chi(\Gamma'_{X,y})-1+\deg\Gamma'_{X,y} \]
In particular, this implies that the coefficient of $\Gamma'\in SS(A)_y$ in the infinite sums (\ref{left action}) and (\ref{right action}) only involves contributions from $\Gamma\in SS(A)$ in a bounded range of degrees and Euler characteristic, and is therefore a finite sum. It follows that the infinite sums in equations  (\ref{left action}) and (\ref{right action}) are well defined.

\begin{lemma}The relations $e^{F}\overleftarrow y=\overrightarrow y e^{F}$, recursively determine $e^{F}$.
\end{lemma}

\pf To start our recursion, we know that $n_{\emptyset}=1$. Choose a $\Gamma_0\in SS(A)$, and assume for induction that we know $n_{\Gamma}$ for all $\Gamma\in SS(A)$ with $\deg\Gamma<\deg\Gamma_0$, or $\deg\Gamma=\deg\Gamma_0$ and $\chi(\Gamma)>\chi(\Gamma_0)$.

Remark \ref{chi bound} implies that, if $n_{\Gamma_0}\neq 0$,  at least one vector in $\Gamma_0$ must be in the form $(a,b)$ with $b<a$. Replace this edge $(a,b)$ from $\Gamma_0$ with an incoming edge $(-a,-1-b)$ to obtain $\Gamma_0'\in SS(A)_{(-a,-1-b)}$. Note that \[\deg \Gamma'_0=\deg\Gamma_0 \text{ and }\chi(\Gamma'_0)=\chi(\Gamma_0)\ .\]
 We now use that the coefficient of $\Gamma_0'$ in $\overrightarrow {(-a,-1-b)}e^{F}-e^{F}\overleftarrow {(-a,-1-b)}$ is 0. Consider the $\Gamma\in SS(A)$ that contribute to this coefficient. Using (\ref{degchi}) we deduce the following.  
 \begin{itemize}
 \item If $\Gamma_{X,(-a,-1-b)}=\Gamma_0'$, then $\deg \Gamma=\deg \Gamma_0$ and $\chi(\Gamma)>\chi(\Gamma_0)$, so we already know $n_\Gamma$.
\item If $\Gamma'_{X,(-a,-1-b)}=\Gamma_0'$, and $(-a,-1-b)+\sum_{v\in X}v=(-1,0)$, then $\deg\Gamma<\deg\Gamma_0$, so we already know $n_\Gamma$.
\item If $\Gamma'_{X,(-a,-1-b)}=\Gamma_0'$,  $(-a,-1-b)+\sum_{v\in X}v=(0,-1)$, and $\abs X>1$, then $\deg \Gamma=\deg \Gamma_0$ and $\chi(\Gamma)>\chi(\Gamma_0)$, so we already know $n_\Gamma$.
\item The only remaining case is when $\Gamma'_{X,(-a,-1-b)}=\Gamma_0'$, and $X$ contains a single vector, equal to $(a,b)$. It follows that $\Gamma=\Gamma_0$.  
\end{itemize}

Therefore, the coefficient of $\Gamma_0'$ in $\overrightarrow {(-a,-1-b)}e^{F}-e^{F}\overleftarrow {(-a,-1-b)}$ consists of known terms plus  the contribution of $\Gamma_0$. This contribution of $\Gamma_0$ is easily computed in terms of $n_{\Gamma_0}$: In this case, $X$ consists of a single vector, equal to $(a,b)$, the number of curves  $n_{X\overleftarrow{(-a,-1-b)}}=a$, and the number of choices of such an $X$ is $\abs{\Aut\Gamma_0}/\abs{\Aut\Gamma_0'}$. It follows that the contribution of $\Gamma_0$ to the coefficient of $\Gamma_0'$ in  $\overrightarrow {(-a,-1-b)}e^{F}-e^{F}\overleftarrow {(-a,-1-b)}$ is $-an_{\Gamma_0}/\abs{\Aut\Gamma_0'}$. As this coefficient is $0$, and all other contributions are known by our inductive hypothesis, this determines $n_{\Gamma_0}$.

The equations $e^{F}\overleftarrow y=\overrightarrow y e^{F}$ therefore determine $e^F$ inductively, as required.

\stop

\subsubsection{Tropical curves for the recursive calculation}\label{tropical description}

\

In this section, we explain how the contact data $\Gamma_{X,y}$ and $\Gamma'_{X,y}$ appearing in our formula for $\overleftarrow y$ and $\overrightarrow y$ correspond to tropical curves. We will use Lemma \ref{gammac} below, but postpone the proof of this Lemma until Section \ref{gf}.

\begin{lemma}\label{gammac}Rigidly constrained tropical curves $\gamma$ in $\expl (M,N)$ with one end constrained and $n_{\gamma}\neq 0$ satisfy the following conditions:
\begin{itemize}
\item  The tropical balancing condition holds at all vertices of $\gamma$ in the interior of $\totb{\expl (M,N)}$.
\item For any vertex $v$ on the positive span of $(1,1)$, the sum of the derivatives of edges leaving $v$ is $(1,0)$.
\item  For any vertex $v$ on the positive span of $(1,1-n)$ (the lower boundary of $\totb{\expl(M,N)}$), the sum of the derivatives of edges leaving $v$ is $(0,1)$.
\item For each vertex $v$ of $\gamma$ in the interior of $\totb{\expl (M,N)}$, let $d_v$ be the valance of $v$ minus $1$, for each vertex on the positive span of $(1,1)$ or  $(1,1-n)$, let $d_v$ be the valence of $v$, and for each vertex $v$ of $\gamma$ sent to $0$, let $d_v=0$. Then, 
\[\sum d_v-\abs{\text{ied}\gamma}\leq 1\]
where $\abs{\text{ied}\gamma}$ is the number of internal edges of $\gamma$.
\end{itemize} 
\end{lemma}

\begin{lemma}\label{gammaXy}Let  $\Gamma=\{\Gamma_1,\dotsc,\Gamma_k\}$ be a rigid tropical curve in $\totb{\expl (M,N)}$ with connected components $\Gamma_i$, and suppose that either
\begin{itemize}
\item $\Gamma_{X,y}$ describes contact data for curves in $\expl(M,N)$,  and $n_{\overrightarrow yX}$ or  $n_{X\overleftarrow y}$ is nonzero; or
\item  $\Gamma'_{X,y}$ describes contact data for curves in $\expl(M,N)$, and either  \begin{itemize}\item  $n_{X\overleftarrow y}$ is nonzero and $y+\sum_{v\in X}v=(-1,0)$, 
\item or $n_{\overrightarrow y X}$ is nonzero and $y+\sum_{v\in X}v=(0,-1)$.\end{itemize}
\end{itemize} Then there exists a unique rigidly constrained tropical curve $\gamma$ in $\totb{\expl (M,N)}$ satisfying the conclusion of Lemma \ref{gammac} with constrained contact data $\Gamma_{X,y}$ or $\Gamma'_{X,y}$ respectively so that
\begin{enumerate}
\item\label{c1} the end labeled by $y$ is constrained above the ray spanned by $-y$ if $n_{X\overleftarrow y}\neq 0$, and is constrained below the ray spanned by $-y$ if $n_{\overrightarrow yX}\neq0$; 
\item\label{c2} there are $k$ vertices $v_i$ sent to $0$ in $\totb{\expl(M,N)}$, and $\gamma_{v_i}=\Gamma_i$;
\item\label{c3} the edges leaving these vertices at $0$ are ends of $\gamma$ if and only if they are not in $X$.

\end{enumerate}
Moreover, besides the vertices sent to $0$, this unique curve $\gamma$ has a unique vertex on each ray spanned by vectors in $X$, and at most one other vertex, which is on the boundary of $\totb{\expl( M,N)}$.  The edges attached to this boundary vertex consist of any edges from $X$ traveling along that boundary, and exactly one other edge (with derivative determined by the balancing condition at this vertex from Lemma \ref{gammac}.) 
The edges attached to any non-boundary vertex on the ray spanned by $w\in X$ consist of 
\begin{itemize}
\item the edges in $X$ traveling along this ray from $0$; 
\item and, one incoming edge and one outgoing edge with incoming and outgoing derivatives 
\[y+\sum_{v\in X\vert v\wedge w>0}v \ \ \ \text{ and }\ \ \ y+\sum_{v\in X\vert v\wedge w\geq 0}v\]
respectively if $n_{\overrightarrow yX}\neq 0$, and incoming and outgoing derivatives
\[y+\sum_{v\in X\vert w\wedge v>0}v\ \ \  \text{ and }\ \ \  y+\sum_{v\in X\vert w\wedge v\geq 0}v\]
if $n_{X\overleftarrow y}\neq 0$.
\end{itemize}

 \end{lemma}

\pf 

We can construct $\gamma$ starting with the constrained end. This travels in from infinity with derivative $y$ until it first hits a  ray spanned by vectors in $X$, where we place a vertex joined to the edges from $X$ traveling along that ray, and with one extra edge leaving with derivative determined by the tropical balancing condition. This leaving edge then travels until it hits the next ray spanned by vectors in $X$, and again interacts with edges from $X$ and leaves with a single edge with derivative determined by the tropical balancing condition. This continues until we have interacted with all the non-boundary edges from $X$, when the final edge leaving either travels to infinity, with derivative $y+\sum_{v\in X}v$, as determined by the tropical balancing condition, or travels to the boundary, where it it is attached to a vertex also attached to any edges from $X$ traveling along that boundary. In this case, $\gamma$ has contact data $\Gamma'_{X,y}$, and our assumptions ensure that this boundary vertex satisfies the balancing condition from Lemma \ref{gammac}.  The  balancing condition ensures that our described edge is always traveling respectively anticlockwise or clockwise, and the condition respectively that $n_{\overrightarrow yX}\neq 0$ or $n_{X\overleftarrow y}\neq 0$,  ensures that our described edge will always hit the next ray from $X$. It should be clear from the construction that $\gamma$ is rigidly constrained, and has vertices and edges as described in Lemma \ref{gammaXy}.

The fact that $\gamma$ is the unique curve satisfying our conditions follows from the following:

\begin{claim}\label{gamma structure} Suppose that $\gamma$ is a rigidly constrained curve satisfying the conditions of Lemma \ref{gammaXy}, and that the constrained edge is not radial, so it is not constrained on a ray from $0\in\totb{\expl(M,N)}$. Then, the non-radial edges of $\gamma$ are connected in a continuous path from the constrained edge to a vertex on the boundary of $\totb{\expl (M,N)}$, or to an end of $\gamma$. Moreover, each ray from $0\in \totb{\expl(M,N)}$ intersects at most one non-radial edge of $\gamma$, and  all vertices of $\gamma$ are either sent to $0$ or are connected to a non-radial edge.  
\end{claim}

First, note that no tropical curve  in $\totb{\expl (M,N)}$ can be rigidly constrained if it has a vertex not sent to zero and only attached to radial edges.  We can deform a tropical curve  in $\expl(M,N)$ continuously by independently scaling each connected component of the curve minus vertices. It follows that for $\gamma$ to be rigidly constrained with the given contact data, when we remove its vertices sent to $0$, it must have a unique connected component  that is not a ray, and this nontrivial component must include the constrained edge.

Starting from any non-radial edge of this nontrival component, the balancing condition implies that we can follow $\gamma$ in a continuous path clockwise through non-radial edges until we arrive at an end or a vertex on the boundary. Similarly, we can follow  anticlockwise until we arrive at an end or boundary vertex. This nontrivial component therefore has at least two ends or vertices on the boundary component. The last condition of Lemma \ref{gammac} then implies that this nontrivial component  must have genus $0$ and exactly two ends or vertices on the boundary.

  As our nontrivial component has zero genus and only two possible ends or boundary vertices to choose from, it follows that all non-radial edges of our nontrivial component are connected in a path traveling either clockwise or anticlockwise from the constrained end to either the boundary vertex or the other end of our nontrivial component. Claim \ref{gamma structure} follows, completing the proof of Lemma \ref{gammaXy}

\stop

\begin{lemma}\label{gamma type} Let $\gamma$ be a rigidly constrained tropical curve in $\totb{\expl(M,N)}$ satisfying the conclusion of Lemma \ref{gammac} and with its constrained edge not on a ray from $0$. Then, there is a unique $\Gamma$, $X$ and $y$ so that  $\gamma$  is the tropical curve constructed in Lemma \ref{gammaXy}.
\end{lemma}

\pf  The rigid tropical curve $\Gamma$ has connected components $\gamma_{v_i}$, where $v_i$ are the vertices of $\gamma$ sent to $0$, and $X$ is the set of ends of $\Gamma$ corresponding to edges that aren't ends, and $y$ is the  incoming derivative of the constrained end of $\gamma$. Claim \ref{gamma structure} implies that all other vertices of $\gamma$ are on a connected path of non-radial edges starting with the constrained edge. These vertices must be connected to vertices at $0$ via the edges in $X$. The tropical balancing condition then implies that $\gamma$ has constrained contact data $\Gamma_{X,y}$ or $\Gamma'_{X,y}$. In particular,  this path of non-radial edges must terminate at an end with derivative $y+\sum_{v\in X}v$, or at a boundary vertex, in which case $y+\sum_{v\in X}v$ must be $(-1,0)$ or $(0,-1)$, in the respective cases when our path is traveling anticlockwise or clockwise, corresponding to the constrained edge being constrained above or below the ray spanned by $-y$. Moreover, the balancing condition ensures that the derivatives of our non-radial edges are as described in Lemma \ref{gammaXy}. For our non-radial edges to reach between the vertices on rays spanned by $X$ with the derivatives described, it is necessary that $n_{\overrightarrow y X}\neq 0$ or $n_{X\overleftarrow y}\neq 0$ in the respective cases when our path is traveling anticlockwise or clockwise. It follows that $\Gamma$, $X$, and $y$ satisfy the requirements of Lemma \ref{gammaXy} and that $\gamma$ is the unique rigidly constrained tropical curve constructed from $\Gamma$, $X$ and $y$ in Lemma \ref{gammaXy}.

\stop

\begin{example}[Computation of some terms in $e^F$]\end{example} 

\begin{figure}[h]
\includegraphics{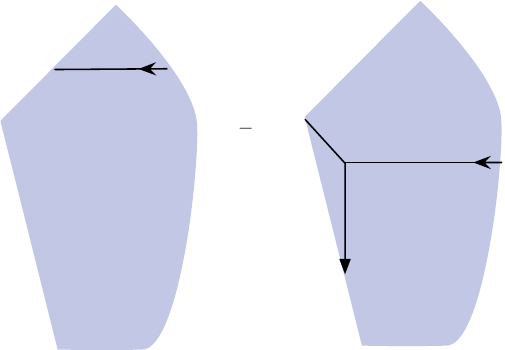}
\caption{Tropical curves contributing to $n_{(-1,0)}$}\label{tec10}
\end{figure}

 To begin with, all we know is that $e^{F}=1\emptyset+\dotsb$. Nevertheless, we can use this information to  compute $n_{(-1,0)}$ using $e^{F}\overleftarrow{(-1,0)}$. The only contributing term is $n_{\emptyset}\emptyset\overleftarrow{(-1,0)}$ --- as pictured in Figure \ref{tec10} on the left --- so $n_{(-1,0)}=1$. Computing using $\overrightarrow{(-1,0)}e^{F}$, the only contributing term is 
\[n_{(1,-1)} \overrightarrow{(-1,0)}(1,-1)= n_{(1,-1)} \{(-1,0)\} \]
as pictured in Figure \ref{tec10} on the right, so $n_{(1,-1)}=n_{(-1,0)}=1$, and \[e^{F}=1\emptyset+1\{(1,-1)\}+\dotsb\ .\]

\begin{figure}[h]
\includegraphics{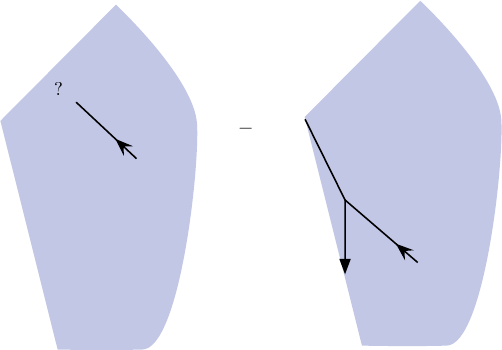}
\caption{Tropical curves contributing to $n_{(-1,g)}$}\label{tec23}
\end{figure}

In fact, $n_{(1,-1)}$ is the only nonzero degree $1$ invariant involved in $e^{F}$. To see that $n_{(1,-1-g)}=0$ for all $g>0$, note that computing $n_{(-1,g)}$ using $\overrightarrow{(-1,g)}e^{F}$ gives $n_{(-1,g)}=n_{(1,-1-g)}$ --- as pictured in Figure \ref{tec23} --- but computing using $e^{F}\overleftarrow{(-1,g)}$ gives  $n_{(-1,g)}=0$, because $(-1,g)$ will never be able to be turned into $(-1,0)$ by interacting with vectors $(a,b)$ with $a>0$, as indicated on the left in Figure \ref{tec23}. 

\begin{figure}[h]
\includegraphics{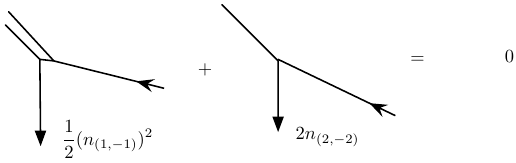}
\caption{Tropical curves contributing to $n_{(-2,1)}$}\label{tec11}
\end{figure}

Similarly, $n_{(-2,1)}=0$ because $(-2,1)$ can't be turned into $(-1,0)$ by interacting with vectors $(a,b)$ so that $a>0$ and $(a,b)\wedge(-2,1)>0$. Calculating $n_{(-2,1)}$ the other way gives  $n_{-2,1}=(n_{{(1,-2)}})^{2}/2 +2n_{(2,-2)}$, so $n_{(2,-2)}=-\frac 14$, as pictured in Figure \ref{tec11}.

In fact, the only nonzero terms in $e^{F}$ with a single vector are $n_{(k,-k)}=\frac{(-1)^{k+1}}{k^{2}}$. This may be verified by calculating $n_{(-k,b)}=0$ for all $b> k$. A nice implication\footnote{See Proposition \ref{abmult} and equation (\ref{dmult}) later in the paper.} of this  is that if $(a,b)$ and $(a+k,b-k)$ are appropriate incoming and outgoing vectors respectively, then 
\[n_{\{(a,b),(a+k,b-k)\}}={\abs{(a,b)\wedge (k,-k)}\choose k}\ .\]

\begin{figure}[h]
\includegraphics{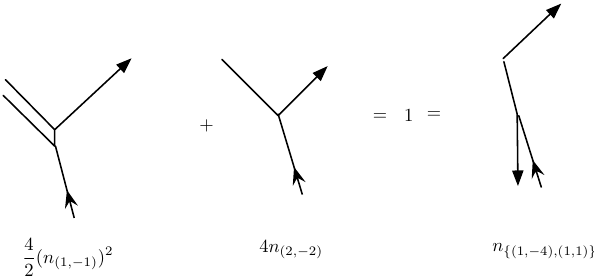}
\caption{Tropical curves contributing to $n_{\{(-1,3),(1,1)\}}$}\label{tec12}
\end{figure}

In particular, $n_{\{(-1,3),(1,1)\}}=1$, so we can calculate $n_{\{(1,-4),(1,1)\}}=1$. The computations involved are pictured in Figure \ref{tec12}.

\section{How recursive calculation of the relative invariants follows from a gluing formula.}\label{gf}

This section explains how equation (\ref{free edge}) follows from our simplified gluing formula, Theorem \ref{simplified gluing}. Applying the explosion functor from \cite{iec}
to the $(M,N)$ introduced in section \ref{relative} gives an exploded manifold $\ex M'$.  The tropical part $\totb{\ex M'}$  of $\ex M'$ is the nonnegative span of $(1,1)$ and $(1,1-n)$ subdivided by the rays $(1,1-k)$.  For us, this subdivison is extraneous information which we can ignore, because $\ex M'$ is a refinement of an exploded manifold $\ex M$,\footnote{We can construct $\ex M$ as the explosion of the toric variety with fan spanned by $(1,1)$, $(1,1-n)$, $(-1,0)$ and $(0,-1)$ relative to the components of the toric boundary divisor corresponding to $(1,1)$ and $(1,1-n)$ .} with tropical part the nonnegative span of $(1,1)$ and $(1,1-n)$, and  Gromov--Witten invariants do not change under the operation of refinement. Accordingly, we shall describe curves in $\ex M$, and draw tropical curves in $\totb{\ex M}$.



For understanding equation (\ref{free edge}), we  regard $n_{\Gamma}$ as a Gromov--Witten invariant of $\ex M$. The curves in $\ex M$ relevant to $n_{\Gamma}$ have tropical part  continuously deformable to  $\Gamma$. In particular, the 
ends  of these tropical curves have derivatives equal  to the vectors in $\Gamma$, and for the $n_{\Gamma}$ from equation (\ref{free edge}), there is a distinguished end with an  extra constraint on its position, and with derivative  $-y$, or `incoming' derivative $y$.  As we shall see, the computation of $n_{\Gamma}$ using $e^{F}\overleftarrow y $ corresponds to constraining this incoming edge  above the ray $-y$, and using $\overrightarrow ye^{F}$ corresponds to constraining this incoming edge  below the ray $-y$.

After choosing where to constrain this incoming edge, the Gromov--Witten invariant $n_{\Gamma}$ decomposes into contributions from rigidly constrained tropical curves. If the incoming edge is constrained to lie on the ray $-y$, then the only such rigid curve is $\Gamma$ itself. Otherwise, there are many possible curves $\gamma$,  each telling how to combine local Gromov--Witten invariants corresponding to vertices to obtain a contribution to $n_{\Gamma}$.  For each vertex $v$ of $\gamma$, the local Gromov--Witten invariant is a Gromov--Witten invariant of the tropical completion $\ex M\tc v$ of $\ex M$ at the vertex $v$, with contact data $\gamma_v$, the tropical curve with a single vertex $v$, and ends with derivatives the derivatives of edges of $\gamma$ leaving $v$. Moreover, we will only need to know the virtual number $n_{\gamma_v}$ of curves in $\ex M\tc v$ with contact data $\gamma_v$ and some collection of ends constrained, as in Definition \ref{constrained count}.

Let us examine curves in $\ex M\tc v$ for the different possible locations of $v$.

\

\noindent{\bf Vertices in the interior of $\totb{\ex M}$}

If $v$ is sent to the interior of $\totb{\ex M}$, \[\ex M\tc v=\ex T^2 \ ,\] because it is $2$-complex-dimensional and has tropical part $\mathbb R^2$. Alternately, we could refine $\ex T^2$ to be the explosion of any compact $2$-complex-dimensional toric manifold. The tropical part of any curve in $\ex T^2$ obeys the usual tropical balancing condition: the sum of derivatives of edges leaving a vertex is $0$.

\begin{lemma}\label{tropical dimension}The virtual fundamental class of the moduli space of curves with genus $g$ and contact data $\gamma_v$ in $\ex T^2$ has real dimension $2(g-1)$ plus twice the valence of $v$. The count of curves with constrained ends and contact data $\gamma_v$ in $\ex T^2$ vanishes unless all but one end is constrained. All such curves have genus $0$.
\end{lemma}

\pf
Because the tangent space of $\ex T^2$ is trivial, the first Chern class of $\ex T^2$ vanishes, so the virtual fundamental class of the moduli space of curves with contact data $\gamma_v$ has real dimension $2(g-1)$ plus twice the valence of $v$.

 If the genus $g>1$, the  dimension of the virtual fundamental class is too high for any constrained count, so it remains to exclude the case $g=1$, which has the correct virtual dimension for a count of curves with all ends constrained. (For our definition of counts of curves with constrained ends, we are assuming all constrained ends have nonzero derivative; see Definition \ref{constrained count}.) To compute the number of such curves, constrain the ends generically. Every balanced tropical curve in $\mathbb R^{2}$  has its infinite edges obeying a one-dimensional constraint caused by the balancing condition. Therefore, there is no balanced tropical curve with all ends constrained generically, so the count of constrained curves vanishes unless $g=0$.
  
 \stop
 \begin{lemma}\label{tropical count} Suppose that $\gamma_v$ has one unconstrained end, and constrained ends with derivative $\alpha$, and $\beta_i$ where all $\beta_i$ are parallel. Then 
 \[n_{\gamma_v}=\prod_i\abs{\alpha\wedge \beta_i}\ .\] 
 \end{lemma}
 
 \pf First, consider the case that $\gamma_v$ has $3$ ends. Analogous to the case of $3$--punctured holomorphic spheres in $(\mathbb C^*)^2$, the moduli space of holomorphic curves in $\ex T^2$ with unconstrained contact data $\gamma_v$ is acted on freely and transitively by the natural action of $\ex T^2$ on itself, and the standard fundamental class of this exploded manifold coincides with the virtual fundamental class of this moduli space.  The evaluation map $ev_e$ at an edge $e$ with derivative $\alpha$ is some constant times the quotient of $\ex T^2$ by the free action with weight $\alpha/\abs\alpha$. It follows that the product of evaluation maps at the two edges with derivative $\alpha$ and $\beta$ has degree $\abs{\alpha\wedge\beta}/\abs{\alpha}\abs{\beta}$, so $n_{\gamma_v}=\abs{\alpha\wedge\beta}$ in this trivalent case. In the more general case, we can constrain the edges with derivative $\beta_i$ to lie on separate lines. There is then a unique rigidly constrained genus $0$ tropical curve with the given constrained contact data. Moreover, this tropical curve is trivalent, and Theorem \ref{simplified gluing} and the trivalent case of our formula imply that $n_{\gamma_v}=\prod_i\abs{\alpha\wedge\beta_i}$, as required. 
 
 \stop  
 
\noindent{\bf{Vertices sent to the boundary of $\totb{\ex M}$}}

 If $v$ is sent to a  $1$--dimensional stratum of $\totb{\ex M}$, \[\ex M\tc v=\ex T\times \expl(\mathbb CP^1,\{0\}) \ ,\] or equivalently, we may use a refinement of $\ex M\tc v$ such as the explosion of  $(\mathbb CP^{1})^{2}$ relative to $3$ of its $4$ toric boundary divisors. 

The tropical part of $\ex T\times\expl(\mathbb CP^1,\{0\})$ is  $\mathbb R\times [0,\infty)$. Because $\totb{\ex M}$ is a subset of $\mathbb R^2$, the tropical part of $\ex M\tc v$ is also a half space of $\mathbb R^2$, but the map between these two half spaces is not the identity map. To specify how the tropical curve $\gamma_v$ in the tropical part of $\ex M\tc v$ translates to contact data for a curve in $\ex T\times\expl(\mathbb CP^1,\{0\})$, we need to know the tropical part of the isomorphism $ \ex T\times \expl(\mathbb CP^1,\{0\})\longrightarrow \ex M\tc v$.

\begin{lemma}\label{boundary identification}If $v$ is on the ray spanned by $(1,1)$, there is an isomorphism $ \ex T\times \expl(\mathbb CP^1,\{0\})\longrightarrow \ex M\tc v$ with tropical part the map
\[(a,b)\mapsto v+a(1,1)+b(1,0) \]
and if $v$ is on the ray spanned by $(1,1-n)$,  there is an isomorphism $ \ex T\times \expl(\mathbb CP^1,\{0\})\longrightarrow \ex M\tc v$ with tropical part the map
\[(a,b)\mapsto v+a(1,1-n)+b(0,1)\]
 \end{lemma}
\pf 

Let us argue the case with $v$ on the ray spanned by $(1,1)$.
 Consider the explosion, $\ex B$ of the toric variety with fan spanned by  $\pm (1,1)$, $(1,1-n)$, $\pm (-1,0)$ and $(0,-1)$, relative to the components of its boundary divisor corresponding to $(1,1)$, $(1,0)$ and $(1,1-n)$. Refining, or otherwise modifying  $\ex M$ away from $v$ does not affect $\ex M\tc v$, so $\ex B\tc v=\ex M\tc v$. Similarly,  $\ex B\tc v=\ex B'\tc v$ where $\ex B'$ is the explosion of the toric manifold with fan spanned by $\pm (1,1)$, and $\pm (1,0)$, relative the components of the toric boundary divisor corresponding to $(1,1)$ and $(1,0)$. Moreover, there is an isomorphism of a refinement of $\ex T\times \expl (\mathbb CP^1,{0})$ with  $\ex B'$, and the tropical part of this isomorphism is $(a,b)\mapsto v+a(1,1)+b(1,0)$. Putting all this together, we get an isomorphism of tropical completions,  
  \[\ex T\times \expl(\mathbb CP^1,\{0\})=\lrb{\ex T\times \expl(\mathbb CP^1,\{0\})}\tc{(0,0)}=\ex B'\tc v=\ex B\tc v=\ex M\tc v\]
and this isomorphism has the require tropical part. The case with $v$ on the ray spanned by $(1,1-n)$ is identical after  swapping the roles of $(1,1-n)$ and $(1,1)$, and $(-1,0)$ and $(0,-1)$.

\stop

\begin{lemma}\label{boundary count} Suppose that $\gamma_v$ is a tropical curve in $\mathbb R\times[0,\infty)=\totb{\ex T\times\expl (CP^1,\{0\})}$ with a single vertex $v$, and no internal edges. Then, the sum of all derivatives of the ends of $\gamma_v$ is $(0,m)$, and the virtual fundamental class of curves with unconstrained contact data $\gamma_v$ and genus $g$ has dimension $2(g-1+m+k)$ where $k$ is the valence of $v$.

 The constrained count of curves with contact data $\gamma_v$ vanishes unless all edges are constrained, and the sum of the derivatives of the ends of $\gamma_v$ is $(0,1)$. In particular, $\gamma$ has a unique end with derivative $\alpha$ not tangent to the boundary of $\mathbb R\times[0,\infty)$, and all other ends have derivative $\beta_i$ tangent to the boundary. Moreover we have
\[n_{\gamma_v}=\prod_i\abs{\alpha\wedge\beta_i}\ ,\]
and all curves contributing to the constrained count have genus 0, and under the projection to  $\expl(\mathbb CP^1,0)$, they have degree $1$.
\end{lemma}

\pf Using the projection to $ \ex T$, and the balancing condition for the tropical part of curves in $\ex T$, we get that the sum of the derivatives of ends of $\gamma_v$ is $(0,m)$. Moreover, under the projection to $\expl (\mathbb CP^1,\{0\})$, our curve has degree $m$. In the case $m=0$, all edges leaving $v$ are parallel to the boundary of $\mathbb R\times[0,\infty)$, so $\gamma_v$ can't be made rigid by constraining ends and $n_{\gamma_v=0}$. We are left with the cases where $m\geq1$. The first Chern class of $\ex T$ vanishes, and the first Chern class of $\expl (\mathbb CP^1,\{0\})$ is represented by intersection with $\{\infty\}\subset\mathbb CP^1$, so the integral of the first Chern class of $\ex M\tc v$ over a curve with tropical part $\gamma_v$ is $m$. The real dimension of the virtual fundamental class of curves with contact data $\gamma_v$ and genus $g$ is therefore $2(g-1+m+k)$ where $\gamma_v$ has $k$ ends. Therefore, the only nonzero constrained count is when $g=0$ and $m=1$.

In the case where $\gamma_v$ has a single end with derivative $(0,1)$, the moduli space of curves with unconstrained contact data $\gamma_v$ is $\ex T$, the virtual fundamental class is represented by this same moduli space, and the evaluation map is the identity $\ex T\longrightarrow \ex T$. This implies that the constrained count in this case is $n_{\gamma_v}=1$.

In the remaining case that $\gamma_v$ has several ends with derivatives $\beta_i$ parallel to the boundary of $\mathbb R\times [0,\infty)$, we can constrain these ends to all lie on a single line off the boundary. There is then a unique rigidly constrained tropical curve with constrained contact data $\gamma_v$ satisfying the tropical balancing condition in the interior and the balancing condition above on the boundary. In particular, in the interior, this tropical curve has a vertex with all our constrained ends entering, and one edge exiting with derivative $(0,-1)$ to terminate at a uni-valent vertex on the boundary. Then Theorem \ref{simplified gluing} together with Lemma \ref{tropical count} and the case of a uni-valent vertex imply that 
\[n_{\gamma_v}=\prod_i\abs{\alpha\wedge\beta_i}\]
as required.

\stop

\noindent{\bf Vertices sent to the corner of $\totb{\ex M}$}

 If $v$ is sent to the corner of $\totb{\ex M}$, tropical completion at $v$ does nothing so $\ex M\tc v=\ex M$.

We do not have a tropical balancing condition at such vertices $v$, however the homology class represented by a curve in $\ex M$ can be read off from its tropical part. Recall that the smooth part, $\totl{\ex M}$, of $\ex M$ is a toric variety with fan $(-1,0)$, $(1,1)$, $(1,1-n)$, and $(0,-1)$, so has a canonical map to $\mathbb CP^2$, which has fan $(-1,0)$, $(1,1)$ and $(0,-1)$. 
\begin{defn}Define the degree of a curve $f:\ex C\longrightarrow \ex M$ to be the degree of the curve in $\mathbb CP^2$ defined by the composition 
\[\totl {\ex C}\xrightarrow {\totl f}\totl{\ex M}\longrightarrow \mathbb CP^2\ .\]
\end{defn}

\begin{lemma} \label{deg and vdim} A curve in $\ex M$ with contact data $\Gamma$ has degree
\[\sum_{(a,b)\in \Gamma} a \ .\]
Moreover,  the virtual fundamental class of the moduli space of curves in $\ex M$ with genus $g$ and contact data $\Gamma$ has real dimension 
\[2(g-1)+2\sum_{(a,b)\in \Gamma}(a+b+1)\]
\end{lemma}

\pf The explosion functor applied to the map $\totl{\ex M}\longrightarrow \mathbb CP^2$ induces a map from $\ex M$ to the explosion of $\mathbb CP^2$ relative to the divisor corresponding to $(1,1)$.  The tropical part of this map sends $(a,b)$ to $a$, so the smooth part of a curve with contact data $\Gamma$ is sent to a curve in $\mathbb CP^2$ of degree $\sum_{(a,b)\in \Gamma}a$, as required.

The first Chern class of $\ex M$ is represented by intersection with the image of the unexploded components of the toric boundary divisor, corresponding to $(-1,0)$ and $(0,-1)$ in the fan of the toric variety we exploded to define $\ex M$. The intersection of a curve with the $(-1,0)$--component  is  the degree of the curve. The intersection with the $(0,-1)$--component is  $\sum_{(a,b)\in\Gamma}b$. So, the integral of the first Chern class of $\ex M$ over a curve with contact data $\Gamma$ is $\sum_{(a,b)\in\Gamma} (a+b)$. 

The real virtual dimension of the moduli space of curves in $\ex M$ with genus $g$ and contact data $\Gamma$ is therefore
\[2(g-1)+2k + 2\sum_{(a,b)\in \Gamma}(a+b)=2(g-1)+2\sum_{(a,b)\in \Gamma}(a+b+1)\]
where, $k$ is the number of vectors in $\Gamma$.
\stop

\

We now  prove Lemma \ref{gammac}. Let $\gamma$ be a rigidly constrained tropical curve in $\totb{\ex M}$ or $\totb{\expl(M,N)}$ with one end constrained. For the tropical part of any holomorphic curve in $\ex M$ or $\expl(M,N)$, the tropical balancing condition holds at vertices in the interior of $\totb{\ex M}$. Moreover, for any vertex $v$ in the interior, the virtual fundamental class of curves with contact data $\gamma_v$ has real dimension at least $2d_v$ where $d_v$ is the valence of $v$ minus $1$.  
 
 If $v$ is a vertex on the boundary Lemma \ref{boundary count} and Lemma \ref{boundary identification} imply that the sum of the derivatives of ends of $\gamma_v$ is  $(m,0)$ or $(0,m)$, in the respective cases that  $v$ is on the ray spanned by $(1,1)$ or $(1,1-n)$. Lemma \ref{dimension vanishing} also implies that if $n_\gamma\neq 0$,  $m\leq 1$, but the case $m=0$ is not possible, because then all edges leaving $v$ would be parallel, and $\gamma$ could not be rigidly constrained. So, we have that the balancing condition for vertices at the boundary of $\totb{\ex M}$  from Lemma \ref{gammac} applies to $\gamma$. For these boundary vertices, the minimum real dimension of the virtual fundamental class is $2d_v$, where $d_v$ is the valence of the vertex $v$. The final condition of Lemma \ref{gammac} follows from our calculation of these $d_v$, and Lemma \ref{dimension vanishing}. This concludes the proof of Lemma \ref{gammac}.

\

\

\begin{prop} The recursion from Section \ref{recursive calculation} holds, and 
\[e^F\overleftarrow y=\sum_{\Gamma\in SS(A)_y}\frac {n_\Gamma}{\abs{\Aut \Gamma}}\Gamma=\overrightarrow y e^F\ .\]
\end{prop}

\pf

We first prove the analogous statement restricting to Gromov--Witten invariants and contact data for $\ex M$ or $\expl (M,N)$, where the number $n$ of blowups is enough that the ray spanned by $-y$ is in the interior of $\totb{\ex M}$, so we can constrain this end above or below the ray spanned by $-y$.

First, our decomposition formula, (\ref{decomposition}), implies that 
\begin{equation}\label{ym0}\sum_{\Gamma'}\frac {n_{\Gamma'}}{\abs{\Aut \Gamma'}}\Gamma'=\sum_\gamma \frac{n_\gamma}{\abs{\Aut\gamma}\abs{\Aut[\gamma]}}[\gamma]\end{equation}
where the first sum is over constrained contact data $\Gamma'$ for  $\ex M$ with a single end with derivative $-y$ constrained above the ray spanned by $-y$, and  second sum is over rigidly constrained tropical curves $\gamma$ in $\totb{\ex M}$ with such contact data, and $[\gamma]$ indicates the constrained contact data of $\gamma$.  Recall that we have used $\Aut\gamma$ for the automorphisms of the tropical curve $\gamma$ fixing the ends of $\gamma$, so $\abs{\Aut\gamma}\abs{\Aut[\gamma]}$ is the size of the automorphism group of $\gamma$ as a tropical curve without labeled ends. 

For any such $\gamma$ with  $n_\gamma\neq 0$, Lemma \ref{gamma type} implies that $\gamma$ is in the form of the tropical curve constructed in Lemma \ref{gammaXy}. In particular, tropical completion of $\gamma$ at $0$ gives a rigid tropical curve $\Gamma$, with a set $X$ of edges from internal edges of $\gamma$, and the contact data of $\gamma$ is $\Gamma_{X,y}$ or $\Gamma'_{X,y}$. Moreover, Theorem \ref{simplified gluing} applies to $\gamma$ with the orientation described in Lemma \ref{gammaXy}, and Lemmas \ref{tropical count}, \ref{boundary identification}, and \ref{boundary count} imply that 
\[n_\gamma=n_\Gamma n_{X\overleftarrow y}.\]
Given $\Gamma$, and $\gamma$ related as above, the number of ways of choosing ends $X$ of $\Gamma$ to create $\gamma$ as in Lemma \ref{gammaXy} is $\abs{\Aut\Gamma}/\abs{\Aut \gamma}\abs{\Aut [\gamma]}$.  
So, 
\begin{equation}\label{ym}\sum_\gamma  \frac{n_\gamma}{\abs{\Aut\gamma}\abs{\Aut[\gamma]}}[\gamma]=\sum_{\Gamma, X}\frac{n_\Gamma}{\abs{\Aut \Gamma}}n_{X\overleftarrow y}\Gamma_{X,y} +\sum_{\Gamma,X}\frac{n_\Gamma}{\abs{\Aut \Gamma}}n_{X\overleftarrow y}\Gamma'_{X,y} \end{equation}
where the sums on the righthand side are over rigid tropical curves $\Gamma$ in $\totb{\ex M}$ and subsets $X$ of the edges of $\Gamma$ so that respectively $\Gamma_{X,y}$  describes contact data for curves in $\ex M$, or $\Gamma'_{X,y}$  describes contact data for curves in $\ex M$ and $y+\sum_{v\in X}v=(-1,0)$. Comparing to (\ref{overleftarrowdef}), the righthand side of (\ref{ym}) is $e^{F_n}\overleftarrow y$, where $F_n$ is the Gromov--Witten potential of $\ex M$ or $\expl(M,N)$, as in Definition \ref{Fndef}. Taking a limit\footnote{We are using that that $n_\Gamma$ does not depend on the number of blowups $n$, so long as the number of blowups is sufficient that $n_\Gamma$ is defined. This can be proved directly, or can be proved using the relations $\overrightarrow ye^{F_n}=e^{F_n}\overleftarrow y$ so long as $n$ is large enough that $-y$ is in the interior of the corresponding $\totb{\expl(M,N)}$. } of the number of blowups $n$ in equations (\ref{ym0}) and (\ref{ym}), we get
\[\sum_{\Gamma\in SS(A)_y}\frac {n_\Gamma}{\abs{\Aut \Gamma}}\Gamma= e^F\overleftarrow y\]
as required. The argument for $\overrightarrow ye^F$ is analogous, except we constrain our end below the ray spanned by $-y$.
     
\stop

\section{Reconstructing the absolute Gromov--Witten invariants from the relative invariants}\label{reconstruct}

To relate our relative Gromov--Witten invariants to the absolute Gromov-Witten invariants of the $n$-fold blowup of $\mathbb CP^{2}$, we shall consider a degeneration of this manifold --- we will then use our relative Gromov-Witten invariants as an essential ingredient in a tropical gluing formula reconstructing the Gromov--Witten invariants of our $n$-fold blowup. As a warmup, we will first describe a corresponding degeneration of $\mathbb CP^{2}$.

\begin{figure}[h]
\includegraphics{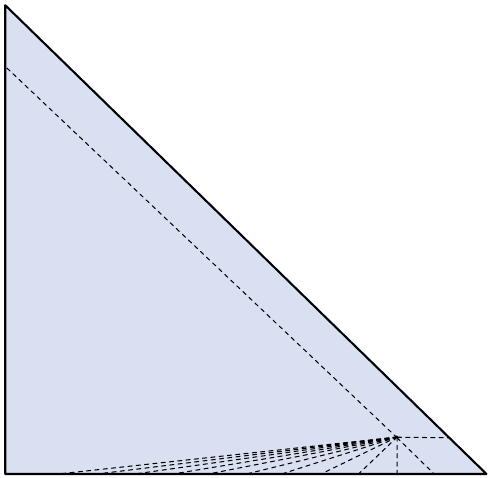}
\caption{Subdivision of the moment polytope of $\mathbb CP^2$}\label{pic1}
\end{figure}

Consider the moment map of $\mathbb CP^{2}$ with the standard torus action. We can  subdivide the image of the moment map as in Figure \ref{pic1} using rays in the directions $(-1,1)$, $(1,0)$, and $(k,-1)$ for all integers $k\in [-n+1,1]$, while ensuring  that all the downward pointing rays intersect the lower edge of the moment polytope.

There exists a toric degeneration of  $\mathbb CP^{2}$ with degenerate fiber consisting of a union of toric manifolds with moment polytopes the pieces of the subdivided triangle in Figure \ref{pic1}, glued along boundary divisors as in this figure. So that we can apply the explosion functor, we want this degeneration to be log smooth when given an appropriate log structure.

\begin{construction}\label{con1}

Consider a dual polytope $P$ to the set of dotted rays in Figure \ref{pic1} --- in other words, consider a convex polygon $P$ with edges orthogonal to the above rays. Choose $P$ to have  integral vertices. 
Then consider the toric partial compactification  $X$ of $(\mathbb C^{*})^{3}$ given 
by the fan consisting of 
\begin{itemize}
\item the cone over $P$, when $P$ is placed in the plane with first coordinate $1$,
\item the cone formed by $(0,-1,0), (0,1,1)$ and the top right vertex of $P$,
\item the cone formed  $(0,1,1)$, $(0,0,-1)$ and the other righthand vertex of $P$,
\item the cone formed by $(0,0,-1)$, $(0,-1,0)$ and top left  vertex of $P$,
\item the cone formed by the righthand face of $P$ and $(0,1,1)$,
\item the cone formed by the top left face of $P$ and $(0,-1,0)$,
\item and the cones formed by $(0,0,-1)$ and all the lower faces of $P$.
\end{itemize}
 Projection of $X$ to the first coordinate  gives the required toric degeneration.
\end{construction}

This toric degeneration $\pi:X\longrightarrow \mathbb C$ is log smooth when $\mathbb C$ is given the log structure from the divisor $0$, and $X$ is given the log structure from the divisor $\pi^{-1}(0)$. This divisor is the union of the toric divisors of $X$ corresponding to all rays in the fan of $X$ apart from $(0,-1,0),(0,1,1), (0,0,-1)$.  To verify that $X$ with this log structure is log smooth, note that the cones formed using two of the directions $(0,-1,0)$, $(0,1,1)$ and $(0,0,-1)$ may be transformed (using an invertible $\mathbb Z$-linear transformation) to standard quadrants, and that each time a cone is formed using a face of $P$ and one of these directions, the configuration formed by the linear plane containing the face of $P$ and the extra direction may be transformed to the standard configuration consisting of the plane spanned by the first two coordinates and $(0,0,1)$.

\begin{construction} \label{con2}

 Blow up $X$ along $n$ complex submanifolds intersecting the singular divisor transversely at $n$ points distributed within the $n$ triangles in Figure \ref{pic1}. By restricting the family $\pi$ to some neighborhood $D$ of $0\in\mathbb C$, we may assume that these $n$ complex submanifolds are transverse to all fibers of $\pi$, so the resulting blown-up family $\pi':X'\longrightarrow D$ is also a log smooth family.

As explained in \cite{iec} or \cite{elc}, we may apply the explosion functor to $\pi':X'\longrightarrow D$ to obtain a smooth family of exploded manifolds.
\[\expl \pi':\expl X'\longrightarrow \expl D\]
\end{construction}

Each exploded manifold has a tropical part consisting of a union of polytopes. The tropical part of $\expl X'$ is the cone over $P$, the tropical part of $\expl D$ is the half line, and the tropical part of $\expl \pi'$ is the projection defining our toric degeneration. The definition of a smooth family in \cite{iec} contains a condition of being surjective on integral vectors --- this condition is satisfied due to our choice of the corners of $P$ having integer coordinates. It is easy to choose a symplectic form taming the complex structure of $\expl X$ in the sense of \cite{cem}. After a choice of symplectic representation of our blowup, this gives a symplectic form on $\expl X'$ taming the complex structure. As the tropical part of $\expl X'$ may be embedded in a quadrant of $\mathbb R^{3}$, the results of \cite{cem} imply that Gromov compactness holds in our family $\expl \pi'$, and we may define Gromov--Witten invariants as in \cite{evc,vfc}.

A general fiber of $\pi'$ is a  $n$-fold blowup of $\mathbb CP^{2}$, and the same is true of the corresponding fiber of $\expl\pi'$. The fiber of $\expl \pi'$ over any point with tropical part $1\in [0,\infty)$ is an exploded manifold $\ex B$ with tropical part $\totb{\ex B}=P$. We shall show how to calculate the Gromov--Witten invariants of $\ex B$. As the Gromov--Witten invariants of exploded manifolds do not change in connected families of exploded manifolds, \cite{evc,vfc}, the Gromov--Witten invariants of $\ex B$ correspond to the Gromov--Witten invariants of $\mathbb CP^{2}$ blown up at $n$ points. 

As explained in \cite{gfgw}, Gromov--Witten invariants of $\ex B$ decompose into a sum of contributions from rigid tropical curves $\gamma$ in the tropical part $\totb{\ex B} $ of $\ex B$. In particular, we can decompose the Gromov--Witten potential from Definition \ref{Gndef} as
\[G_n:=\sum_{g,\beta}n_{g,\beta}x^{g-1}q^\beta=\sum_\gamma \frac{n_\gamma}{\abs{\Aut\gamma}} x^{g(\gamma)-1}q^{\beta(\gamma)}\]
where $n_\gamma$ is the virtual number of rigid curves in $\ex B$ with an isomorphism of their tropical part with $\gamma$, defined as in (\ref{ngammadef}),  and the curves contributing to $n_\gamma$  have genus $g(\gamma)$ and represent homology class $\beta(\gamma)$. We will be able to compute $n_\gamma$ using Theorem \ref{simplified gluing}

\begin{figure}[htbp]
\includegraphics{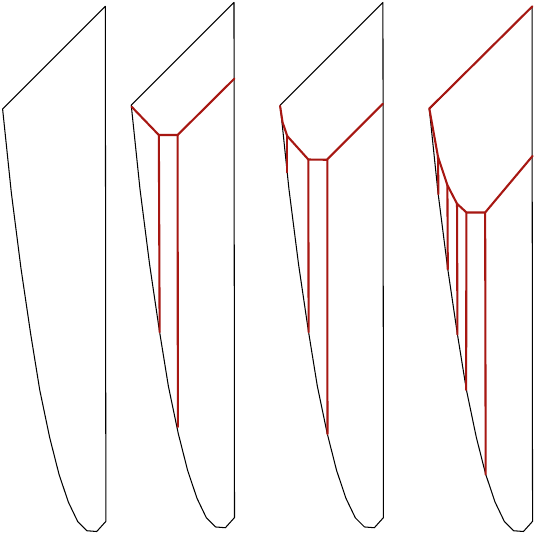}
\caption{Tropical curves in the tropical part of $\ex B$}\label{tec2}
\end{figure}

Figure \ref{tec2} shows $\totb{\ex B}$, and some tropical curves in $\totb{\ex B}$. The left and righthand curves both contribute $1$ to the Gromov--Witten count of curves and the middle picture does not contribute, because it is not rigid --- actually this middle tropical curve deforms to the lefthand tropical curve. 

\begin{remark}Confusingly,  when using the non-generic complex structure on $\ex B$ from Constructions \ref{con1} and \ref{con2}, there are no genuine holomorphic curves with tropical part given by the left and righthand pictures, but there are holomorphic curves  with tropical parts such as those pictured in the middle that deform to the lefthand picture (and something similar happens for the righthand picture).   When a generic almost complex structure on $\ex B$ is used,  there are unique holomorphic  curves in $\ex B$ with tropical parts the left and righthand pictures, but there does not exist any holomorphic curve with tropical part given by the middle picture. 
\end{remark}

Use $H$ to denote the homology class represented by the pullback of a line to the $n$-fold blowup of $\mathbb CP^{2}$, and let $E_{i}$ be the homology class of the $i$th exceptional divisor. As proved in \cite{dre}, our exploded manifold $\ex B$ has the same DeRham cohomology as $\mathbb CP^{2}$ blown up at $n$ points, so the same classes make sense in $\ex B$ --- we can choose a representative for $H$ with tropical part the lefthand corner and top boundary of $\totb{\ex B}$, and a representative for each $E_{i}$ with tropical part the $i$th corner at the bottom of $\totb{\ex B}$. With these choices, we can measure the homology class represented by curves using their intersection with our representatives, and talk of the individual contributions of vertices of tropical curves to the overall homology class. 
From the left, the first and second curves are rational curves in the class $H-E_{3}-E_{5}$, and the last is a rational curve representing  $2H-E_{1}-E_{2}-E_{3}-E_{4}-E_{6}$.

\

As proved in \cite{gfgw}, the Gromov--Witten contribution  of each rigid tropical curve in $\totb{\ex B}$  is determined by taking a fiber product of relative Gromov--Witten invariants corresponding to its vertices. These relative invariants at a vertex $v$ are Gromov--Witten invariants of the tropical completion $\ex B\tc v$ of $\ex B$ at $v$.  We must consider curves in   $\ex B\tc v$ for vertices of the following types:
 
 \begin{description}
 \item[Type A] vertices in the interior of $\totb{\ex B}$;
\item[Type B] vertices in $1$--dimensional strata of $\totb{\ex B}$; 
\item[Type C] vertices on the lefthand corner of $\totb{\ex B}$ in Figure \ref{tec2};
\item[Type D] vertices on one of the two righthand corners of $\totb{\ex B}$ in Figure \ref{tec2};
\item[Type E] vertices on one of the $n$  bottom corners of $\totb{\ex B}$ in Figure \ref{tec2}.
 \end{description}

For $v$ a Type A vertex, $\ex B\tc v=\ex T^2$, so the tropical part of any holomorphic curve in $\ex B$ satisfies the tropical balancing condition at these vertices, and we can analyse the Gromov--Witten invariants of $\ex B\tc v$ using lemmas \ref{tropical dimension} and \ref{tropical count}.

For $v$ a Type B vertex, $\ex B\tc v=\ex T\times\expl(\mathbb CP^1,0)$, and we can analyse Gromov--Witten invariants using Lemma \ref{boundary count}.

\begin{lemma}\label{Bvertex} If $v$ is a Type B vertex, there is an isomorphism $\ex T\times\expl(\mathbb CP^1,0)\longrightarrow \ex B\tc v$, with tropical part sending $(0,b)$ to 
\begin{itemize}
\item $v+b(1,0)$ if $v$ is on the top left boundary of $\totb{\ex B}$ in Figure \ref{tec2};
\item $v+b(-1,-1)$ if $v$ is on the righthand boundary of $\totb{\ex B}$ in Figure \ref{tec2};   
\item $v+b(0,1)$ if $v$ is on one of the remaining $1$--dimensional strata of $\totb{\ex B}$, all of which are on the bottom of $\totb{\ex B}$ in Figure \ref{tec2}.
\end{itemize}
\end{lemma}

\pf Our point $v\in \totb{\ex B}$ gives a point in the tropical part of the total space of the degeneration $\expl \pi'$ from Construction \ref{con2}. Consider taking the tropical completion of the degeneration $\expl \pi'$ at $v$. Because $v$ is in a stratum that avoids all the blowups of Construction \ref{con2}, this tropical completion is equal to the tropical completion of $\expl \pi$ at $v$, where $\pi: X\longrightarrow \mathbb C$ is the toric model degeneration of Construction $\ref{con2}$, and the explosion is relative to the divisor $\pi^{-1}(0)$. Our  $v$ is on a $1$--dimensional stratum of the polytope $P$, corresponding to a $2$--dimensional stratum of $\totb{\expl X}$ and a $2$--dimensional stratum of the toric fan of $X$. 

Let $e\subset P\subset \mathbb R^3$ indicate the closure of the stratum of $P$ containing $v$, and let $\alpha$ be $(1,0)$, $(-1,-1)$ or $(0,1)$ in the respective cases listed above.  Using that any modification of $\expl X$ away from $v$ does not affect $\expl X\tc v$, we get that $\expl X\tc v=\expl Y\tc v$, where $Y$ is a toric variety with fan the two cones spanned by $e$ and $\pm \alpha$, and the explosion is relative to every toric boundary stratum of $Y$ apart from the component corresponding to $-\alpha$. This exploded manifold $\expl Y$ is isomorphic to $\expl Y_0\times \expl(\mathbb CP^1,0)$, where $Y_0$ is some $2$--dimensional toric variety, and the tropical part of this isomorphism sends the cone on $e$ to  $\totb{\expl Y_0}$ , and sends the ray spanned by $\alpha$ to $\totb{\expl (\mathbb CP^1,0)}$. Tropical completion of $\expl Y_0\times \expl(\mathbb CP^1,0)$ at the image of $v$ is $\ex T^2\times \expl(\mathbb CP^1)$, so we obtain an isomorphism $\ex T^2\times \expl (\mathbb CP^1)\longrightarrow \expl X\tc v$ with tropical part sending $(0,0,b)$ to $v+b\alpha$. The tropical completion of $\expl\pi$ at $v$ is then isomorphic to the projection, $\ex T^2\times\expl (\mathbb CP^1)\longrightarrow \ex T$,  onto the first coordinate. Our $\ex B\tc v$ is therefore isomorphic to the restriction of this map to a fiber, which is $\ex T\times\expl(\mathbb CP^1,0)$, and the tropical part of this isomorphism sends $v+b\alpha$ to $(0,b)$, as required.

\stop

\begin{lemma}\label{CDvertex} If $v$ is a Type C or D vertex, there is an isomorphism of $\ex B\tc v$ with the explosion of the toric variety with moment polytope a quadrilateral in Figure \ref{pic1}, relative to its dotted boundary divisors. Moreover, the tropical part of this isomorphism is the translation sending $v$ to $0$. If $v$ is in left, top right, or bottom right corner of \totb{\ex B} in Figure \ref{tec2}, we respectively use the bottom left, top, and right quadrilaterals in Figure \ref{pic1}. In particular, if $v$ is Type C,  $\ex B\tc v$ is the exploded manifold $\ex M$ from Section \ref{gf}, with refinement $\expl (M,N)$ from Section \ref{relative}. 
\end{lemma}

\pf

The smooth part of $\ex B$ is the complex of toric manifolds with moment polytopes shown in Figure \ref{pic1}, blown up once inside the inverse image of each of the triangles. In particular, the smooth part of the stratum with tropical part $v$ is the inverse image of the  corresponding quadrilateral from Figure \ref{pic1} minus its dotted boundry divisors. The smooth part of $\ex B\tc v$ is the toric variety with moment polytope this quadrilateral, with the stratification given by the dotted boundary divisors. It follows that $\ex B\tc v$ is the explosion of this toric variety relative to the dotted components of its toric boundary divisor. 

\stop

\begin{cor}\label{Ddimension} If $v$ is a type $D$ vertex, then the minimal real dimension of the moduli space of curves with unconstrained contact data $\gamma_v$ is $2d_v$, where  
\[d_v=-1+\sum_{(a,b)\in \gamma_v}(1+a-2b)\]
if $v$ is at the top right corner of $\totb{\ex B}$ in Figure \ref{tec2}, and
\[d_v=-1+\sum_{(a,b)\in\gamma_v}(1-b)\]
if $v$ is at the bottom right corner. In the above, $(a,b)\in \gamma_v$ indicates that $(a,b)$ is a derivative of an edge leaving $v$. This minimal dimension is only achieved for the zero genus component of the moduli space of curves with contact data $\gamma_v$. 

In particular, so long as $\gamma_v$ has at least one nonconstant end, $d_v$ is at least the valence of $v$. 
\end{cor}

\pf In each case, $\ex B\tc v$ is the explosion of the toric manifold with moment polytope the relevant quadrilateral in Figure \ref{pic1}, relative to the dotted boundary components. The Chern class of $\ex B\tc v$ is therefore represented by intersection with the solid boundary components. For curves with contact data $\gamma_v$, this intersection is 
\[\sum_{(a,b)\in\gamma_v}(a-2b)\ \ \ \text{ or }\sum_{(a,b)\in\gamma_v}-b\text{ respectively.}\]
The calculation of $d_v$ then follows from the general formula (\ref{vdim}) for the dimension of the virtual fundamental class with genus $0$ --- higher genus gives higher dimension. Note that the ends of $\gamma_v$ have derivatives in a cone where, if they are nonzero,  their contribution to the above sum is strictly positive. It follows that $d_v$ is at least the valence of $v$ if $\gamma_v$ has a nonconstant end. 

\stop 

\begin{lemma}\label{Evertex} If $v$ is a Type E vertex, $\ex B\tc v$ is isomorphic to the explosion of $\mathbb CP^2$ blown up at one point, relative the union of two lines.
\end{lemma}

\pf As argued in the proof of Lemma \ref{CDvertex}, the smooth part of $\ex B\tc v$ is the toric manifold with moment polytope one of the small triangles in Figure \ref{pic1} blown up once somewhere away from the dotted boundary divisors, and with the stratification from the dotted toric boundary divisors. Moreover, $\ex B\tc v$ is the explosion of this manifold relative to these two divisors.  This toric manifold is isomorphic to $\mathbb CP^2$ and each dotted boundary divisor is isomorphic to a line, so $\ex B\tc v$ is isomorphic to the explosion of $\mathbb CP^2$ blown up at one point, relative to the union of two lines. 

\stop

\begin{lemma}\label{Edimension}Suppose that $\gamma$ is the tropical part of a curve in $\ex B$, and $v$ is a Type E vertex of $\gamma$. Then the sum of the derivatives of edges exiting $v$ is $(0,d)$, where $d\geq 0$. Moreover, the real dimension of any nonempty component of the virtual fundamental class of the moduli space of curves in $\ex B\tc v$ with contact data $\gamma_v$ is at least $2d_v$, where $d_v$ is the valence of $v$ minus $1$. The component of the moduli space of curves with contact data $\gamma_v$ achieving this minimum dimension consists of curves of genus $0$ representing a homology class that intersects the exceptional sphere in $\ex B\tc v$ $d$ times.  
\end{lemma}

\pf We have that $\ex B\tc v$ is the explosion of the $1$-point blowup of $\mathbb CP^2$ relative to two lines.   In the tropical part of $\ex B\tc v$, these two lines correspond to the two $1$--dimensional boundary strata of $\ex B\tc v$, which are rays traveling in the directions $(-1,n-l)$ and $(1,1+l-n)$, for $v$ at the $l$th bottom corner of $\totb{\ex B}$. The fact that any homology class in $\totl{\ex B\tc v}$ has equal intersection with our two lines then implies our balancing condition that the sum of the derivatives of edges exiting $v$ is $(0,d)$. Moreover,  these edge derivatives have to be in the positive span of $(-1,n-l)$ and $(1,1+l-n)$, so $d$ is strictly positive unless the derivatives of all exiting edges vanish.

The first Chern class of $\ex B\tc v$ is represented by the Poincare dual of the strict transform of a line passing through the blown up point. In particular the first Chern class evaluated on any holomorphic curve is non-negative. It follows from (\ref{vdim}) that our minimum dimension for the moduli space of curves in $\ex B\tc v$ with contact data $\gamma_v$ is $2d_v$ where $d_v$ is the valence of $v$ minus $1$. This minimum dimension is only achieved when we restrict to curves with genus $0$ on which the first Chern class vanishes. If the sum of the derivatives of edges of $\gamma_v$ is $d$, such curves intersect the exceptional sphere $d$ times --- they represent $d$ times the pullback of a line from $\mathbb CP^2$ minus $d$ times the exceptional sphere. 

\stop

\begin{lemma}\label{ga} Suppose that $\gamma$ is a connected tropical curve in $\totb{\ex B}$. Then $n_\gamma=0$ unless $\gamma$ is a single point of Type E, or all the following conditions hold.
\begin{enumerate}[i]
\item\label{ga1} $\gamma$ has genus $0$.
\item\label{ga2} $\gamma$ has a uniqe Type C vertex.
\item\label{ga3} Each connected component of $\gamma$ minus its Type C vertex has a unique vertex on the closure of the righthand boundary of $\totb{\ex B}$ in Figure \ref{tec2}. 
\item\label{ga4} This unique vertex on the closure of the righthand boundary of $\totb{\ex B}$ is univalent, with exiting derivative $(-1,-1)$.
\item\label{ga5} Each connected component of $\gamma$ minis its Type C vertex has a connected path of edges with derivative $(1,*)$, and  all other edges are vertical and join a univalent vertex of Type E  to a vertex on this path.
\end{enumerate}
\end{lemma}

\pf

First consider the case that $\gamma$ is a single point. For $\gamma$ to be rigid, it must  be at a corner of $
\totb{\ex B}$; Lemma \ref{CDvertex} implies that $n_\gamma=0$ if $\gamma$ is sent to a corner of Type C or D, therefore the remaining case is that $\gamma$ is a single vertex of Type E.

From now on, we assume that $\gamma$ is rigid and has at least one edge. For $\gamma$ to be rigid, no internal edges can be constant. Moreover, if $n_\gamma\neq0$,  $\gamma$ has no constant ends, because the virtual fundamental class with $n$ constant ends is a $n$--complex dimensional bundle over the corresponding virtual fundamental class with no constant ends, and this negative-dimensional virtual fundamental class vanishes. Therefore, $\gamma$ has no ends, and each edge of $\gamma$ has nonzero derivative.

 Now consider the minimal real dimension, $2d_v$ of curves with contact data $\gamma_v$. If $v$ is on a $1$--dimensional stratum, at least one edge of $v$ must leave that stratum for $\gamma$ to be rigid. Therefore, lemmas \ref{Bvertex} and \ref{boundary count} imply that $d_v$ is at least the valence of $v$. Lemmas \ref{tropical dimension} and \ref{Edimension} and Corollary \ref{Ddimension} imply the following bounds on $d_v$:
 \begin{description}
 \item[Type A] $d_v=\text{valence}(v)-1$
 \item[Type B] $d_v\geq\text{valence}(v)$
 \item[Type C] $d_v\geq 0$
 \item[Type D] $d_v\geq\text{valence}(v)$
 \item[Type E] $d_v\geq \text{valence}(v)-1$
 \end{description}

The balancing conditions at vertices of Type A, B and E imply that each edge of $\gamma$ whose derivative has nonzero first component must be connected through a path of such edges to a vertex of Type B or D on the righthand boundary of $\totb{\ex B}$. Moreover, Lemma \ref{dimension vanishing}, Corollary \ref{Ddimension}, and lemmas \ref{Bvertex} and \ref{boundary count} imply that the derivative of the edge hitting the righthand boundary must be $(1,*)$, and that any other edges exiting this boundary vertex must be parallel to the boundary. Our balancing condition at vertices of Type A, B, and E therefore implies that there are at least as many vertices of Type B and D as there are edges leaving vertices of Type C, and that 
\[\abs{\text{Type B or D}}\geq\sum_{v\in \text{Type C}}\sum_{(a,b)\in\gamma_v}a\]
using the obvious notation for the set of vertices of a given type, and the set of derivatives of ends of $\gamma_v$.
Lemma \ref{dimension vanishing} implies that $\sum d_v$ is less than or equal to the number, $\abs{\text{ied}\gamma}$, of edges of $\gamma$. 

We have that 
\begin{equation}\label{dvcount}\begin{split}\sum_vd_v\geq  &\sum_v(\text{valence} (v)-1)+\abs{\text{Type B or D}}-\sum_{v\in \text{Type C}}(\text{valence}(v)-1)
\\ &= \abs{\text{ied}\gamma}+ \text{genus}(\gamma) -1 +\abs{\text{Type B or D}}+\abs{\text{Type C}}-\sum_{v\in \text{Type C}}\sum_{(a,b)\in\gamma_v}1  \end{split}\end{equation}
where $\text{genus}(\gamma)$ is the genus of the graph $\gamma$, not to be confused with the genus,  $g(\gamma)$, of the curves counted by $n_\gamma$.

Our balancing conditions imply that if there are no Type C vertices, there are at least two Type B or D vertices, in which case Lemma \ref{dimension vanishing} and (\ref{dvcount}) imply that $n_\gamma=0$. Therefore, there is at least one Type C vertex, and to get $\sum_v d_v\leq \abs{\text{ied}\gamma}$, we require that 
\begin{itemize}
\item $\text{genus}(\gamma)=0$;
\item there is a unique Type C vertex;
\item each edge exiting this Type C vertex has derivative $(1,*)$, and is connected to a unique Type B or D vertex on the righthand boundary of $\totb{\ex B}$ through a path of edges not passing through the Type C vertex;
\item and there are no other Type B or D vertices.
\end{itemize}
In particular, we have verified items \ref{ga1}, \ref{ga2} and \ref{ga3}. The above, and our observations about $d_v$ for vertices on the righthand boundary of $\ex B$ imply that each vertex on the righthand boundary must be univalent, with entering derivative $(1,1)$, so item \ref{ga4} holds. Moreover, the above and our balancing conditions imply that each connected component of $\gamma$ minus its unique Type C vertex has a connected path of edges with derivative $(1,*)$ which terminates at the vertex on the righthand boundary of $\totb{\ex B}$, and that any other edge other edge must be vertical. As there are no other vertices of Type B or D, the fact that $\gamma$ is rigid implies that all these vertical edges must connect a vertex on our path of horizontal edges to a Type E vertex.  As $\gamma$ has zero genus, these Type E vertices must be univalent. Therefore item \ref{ga5} holds.

\stop

The following is an easy corollary of Lemma \ref{ga} and our analysis of $d_v$ for vertices of different types.
\begin{cor} \label{gamma orientation}If $\gamma$ is a tropical curve in $\totb{\ex B}$ with $n_\gamma\neq0$, then Theorem \ref{simplified gluing} applies to calculate $n_\gamma$, with vertical edges of $\gamma$ oriented up, and all other edges oriented to have derivative $(1,*)$.

\end{cor}

\begin{lemma} \label{gamma homology}Suppose that $\gamma$ is a connected tropical curve in $\totb{\ex B}$ with $n_\gamma\neq 0$. Then the genus and homology class of curves contributing to $n_\gamma$ are determined as follows: 

If $\gamma$ consists of a single point, sent to the $i$th Type E vertex, $n_\gamma=1$, and counts the $i$th exceptional sphere.

Otherwise, let $(1,1-m_j)$ for $j=1,\dotsc,k $ be the list of derivatives of edges leaving the Type C vertex of $\gamma$  and let $(0,d_i)$ be the sum of the derivatives of edges leaving vertices at the $i$th Type E corner of $\totb{\ex B}$. Then curves contributing to $n_\gamma$ have genus 
\[ g(\gamma):= 1+ \sum_{j=1}^k(m_j-3)  \]
and represent the homology class
\[\beta(\gamma):=kH-\sum_id_iE_i\]
where $H$ is the class represented by a pullback of a line from $\mathbb CP^2$, and $E_i$ is the class represented by the $i$th exceptional sphere. 
\end{lemma}

\pf

Let us first consider the case that $\gamma$ is a point mapping to the $i$th Type E corner of $\totb{\ex B}$. This is the tropical part of a curve compactly supported in the corresponding stratum of $\ex B$, which Lemma \ref{Evertex} tells us is isomorphic to  the blowup of $\mathbb CP^2$ at one point, with two lines removed. This curve must therefore represent the homology class of some multiple of the exceptional sphere. The first Chern class evaluated on this exceptional sphere is $1$, so the only way for such a curve to be rigid is to have genus $0$ and represent the homology class of the exceptional sphere. The only such curve is the exceptional sphere itself.

Let us now restrict to the case that $\gamma$ is nonconstant.

As specified in lemmas \ref{tropical count}, \ref{boundary count}, \ref{Edimension}, and Corollary \ref{Ddimension}, the genus of curves contributing to $n_{\gamma_v}$ is zero for all vertices apart from the unique type C vertex. Moreover, Corollary \ref{gamma orientation} gives that at the Type C vertex, we are using unconstrained contact data, so Lemma \ref{deg and vdim} gives that the rigid curves we are counting at this type C vertex have genus 
\[1-\sum_{j=1}^k(1+(1-m_i)+1) =1+\sum_{j=1}^k(m_i-3)\ .\]
As specified by Lemma \ref{ga}, $\gamma$ has zero genus, so the above is the genus of the curves contributing to $n_\gamma$.

Now we consider the homology class represented by our curves. We can choose complex submanifolds of $\ex B$ representing $H$ and $E_i$. In particular, let $\tilde H$ be the component of the toric boundary divisor corresponding to the ray $(0,-1,0)$ in Construction \ref{con1}, and let $\tilde H'$ be its image in the total space of the family $\expl \pi'$ from Construction \ref{con2}. The restriction of $\expl \pi'$ to $\tilde H'$ is still a family, and the intersection of $\tilde H'$ with a generic fiber represents the pullback of a line from $\mathbb CP^2$, so the intersection $H$  of $\tilde H'$ with the fiber $\ex B$ also represents this class. Thinking of Figure \ref{pic1} as a schematic representation of the smooth part of $\ex B$,  this $H\subset \ex B$ is the complex submanifold with image in the smooth part of $\ex B$ the lefthand boundary divisor in Figure \ref{pic1}.   A curve with tropical part $\gamma$ can potential intersect $H$ at any vertex on the top left boundary of $\totb{\ex B}$. As there are no vertices on the interior of this stratum, and any vertex on the top corner of $\totb{\ex B}$ is univalent with incoming derivative $(1,1)$ our curve only intersects $H$ at its Type C vertex. The intersection with $H$ here is computed in Lemma \ref{deg and vdim} as the sum of the fist components of derivatives of edges leaving our Type C vertex. So, our curve's intersection with $H$ is  $k$.

Similarly, we let $\tilde E_i$ be the ith exceptional divisor from Construction \ref{con2}. Again, $\expl \pi'$ restricted to $\tilde E_i$ is still a family. The intersection of $\tilde E_i$ with a generic fiber is the $i$th  exceptional sphere, and the intersection of $\tilde E_i$ is a complex submanifold $E_i$ of $\ex B$ supported at the $i$th Type E corner of $\totb{\ex B}$, and equal to the exceptional sphere under the identification of this stratum with the blowup of $\mathbb CP^2$ minus $2$ lines. Lemma \ref{Edimension} then implies that curves contributing to $n_\gamma$ intersect $E_i$ $d_i$ times, where $(0,d_i)$ is the sum of the derivatives of edges leaving the $i$th Type E corner of $\totb{\ex B}$. 

Intersection with $H$ and $E_i$ determine second homology in blowups of $\mathbb CP^2$; our curves contributing to $n_\gamma$ represent the homology class
\[kH-\sum_id_iE_i\]
as required.

\stop

To determine $n_\gamma$, we still require $n_{\gamma_v}$ for Type E vertices $v$. This is the subject of Section \ref{last}. We need these invariants for univalent $v$ with outgoing derivative $(0,d)$, and in Section \ref{last}, we calculate these invariants as
\begin{equation}\label{n0d}n_{(0,d)}=\frac{(-1)^{d+1}}{d^{2}}\ .\end{equation}
Equation (\ref{n0d}) follows from lemmas \ref{dmult} and \ref{AA'}. 

\begin{prop}\label{GnFn proof} The equation (\ref{GnFn}), relating our absolute and relative Gromov--Witten invariants, holds. Moreover, in the equation
\[G_n:=\sum_{g,\beta}n_{g,\beta}x^{g-1}q^\beta=\sum_\gamma \frac{n_\gamma}{\abs{\Aut\gamma}} x^{g(\gamma)-1}q^{\beta(\gamma)}\]
we may take the sum over $\gamma$ satisfying the conditions of Lemma \ref{ga} so that if any vertex $v$ of Type A or E has an edge with incoming derivative $(1,k)$, it has an edge with outgoing derivative $(1,k+1)$. Moreover, for such $\gamma$, 
\[n_\gamma=n_{\gamma_{v_0}}\]
where $v_0$ is the unique Type C vertex of $\gamma$. 
\end{prop}

\pf

First, suppose that a tropical curve satisfies the conditions of Lemma \ref{ga}, but has some Type A or E vertex $v$ with incoming derivative $(1,k)$ and and outgoing derivative $(1,k+d)$, with $d>1$. If $v$ is a Type A vertex, there are several possibilities for how $v$ is attached to Type E vertices, labeled by partitions $\lambda$ of $d$ so that the vertical edges attached to $v$ have  incoming derivatives $(0,\lambda_i)$. We will show that the contribution of these different tropical curves cancel. In particular, for each partition $\lambda$ of $d$, let $\gamma^\lambda$  be a tropical curve, with each $\gamma^\lambda$ identical, except for the vertical edges attached to $v$, which have derivatives $(0,\lambda_i)$. Furthermore, define $\Aut'\gamma^\lambda$ to be the group of automorphisms of $\gamma^\lambda$ that only permute these vertical edges.
\begin{claim}\label{reduced cancelation}\[\sum_\lambda \frac{n_{\gamma^\lambda}}{\abs{\Aut'\gamma^\lambda}}=0\]
\end{claim}

At the Type E vertex  attached to the edge with derivative $(0,\lambda_i)$, equation (\ref{n0d}) gives a local invariant 
$(-1)^{\lambda_i+1}/\lambda_i^2$. At $v$, Lemma \ref{tropical count} gives a local invariant of $\prod_i\lambda_i$. Theorem \ref{simplified gluing} then gives 
\[n_\gamma^\lambda=c(-1)^{\sum_i 1+\lambda_i}\prod_i\frac 1{\lambda_i}=c\sigma(\lambda)\prod_i\frac 1{\lambda_i}\]
where $\sigma(\lambda)$ is the sign of the permutation of $d$ with cycle type $\lambda$ and $c$ is the product of $n_{\gamma^\lambda_v}$ for vertices $v$ not attached to our vertical edges --- so $c$ is independent of $\lambda$. When the partition $\lambda$ has $n_k$ elements of size $k$,  we have  $\abs{\Aut'\gamma^\lambda}=\prod_k n_k!$.  The stabiliser group of a permutation of cycle type $\lambda$ has size $z_\lambda=\abs{\Aut'\gamma^\lambda}\times \prod_i\lambda_i$. For $d>1$, $\sigma$ is a surjective homomorphism to $\mathbb Z_2$ so the number of even permutations equals the number of odd permutations. This implies  the  identity
\[ \sum_\lambda \sigma(\lambda)/z_\lambda=0 \]
so
\[\sum_\lambda \frac{n_{\gamma^\lambda}}{\abs{\Aut'\gamma^\lambda}}=\sum_\lambda c\frac{\sigma(\lambda)}{z_\lambda}=0 \ .\]

With Claim \ref{reduced cancelation} now proved,  define $\Aut_0\gamma$ to be the group of automorphisms of $\gamma$ fixing edges leaving the unique Type C vertex. Group the contributions of tropical curves according to their contact data at their Type C vertex. We have
\begin{equation}\label{GnFn1}G_n=\sum_{i=1}^n x^{-1}q^{E_i}+\sum_\Gamma \frac 1{\abs{\Aut \Gamma}}\sum_{\gamma\vert \gamma_{v_0}=\Gamma}\frac{n_\gamma}{\abs{\Aut_0\gamma}}x^{g(\gamma)-1}q^{\beta(\gamma)}\end{equation}
where the righthand sums are over unconstrained contact data $\Gamma$ at the Type C vertex, and tropical curves $\gamma$ with a given identification of $\gamma_{v_0}$ with $\Gamma$. Because the automorphisms in $\Aut_0\gamma$ just permute vertical edges attached to the same Type A vertex, $\Aut_0\gamma^\lambda/\Aut'\gamma^\lambda$ is independent of $\lambda$, so the contributions of $\gamma^\lambda$ to the above sum cancel.

We have shown the cancelation of the contributions of tropical curves with Type A vertices having incoming and outgoing edges with derivatives $(1,k)$ and $(1,k+d)$ for $d>1$.  Similarly, we have the following: 
\begin{claim}\label{E reduction} if $\gamma$ has a Type E vertex $v$ with incoming and outgoing  edges with derivative $(1,k)$ and $(1,k+d)$, where $d>1$ then $n_{\gamma_v}=0$ so  $n_\gamma=0$.
\end{claim}

To prove Claim \ref{E reduction}, compute $n_{\gamma_v}$ by constraining the incoming edge above the ray emanating from $v$. Then the tropical curves contributing to this count have one type $A$ vertex joined to vertical edges from $v$ with derivatives $(0,\lambda_i)$ from some partition $\lambda$ of $d$. Moreover, Theorem \ref{simplified gluing}, Lemma \ref{tropical count}, and (\ref{n0d}) give
\[n_{\gamma_v}=\sum_\lambda\frac{\sigma(\lambda)}{z_\lambda}=0\]
as required.

We are left with the contributions of tropical curves $\gamma$ such that every vertical edge has derivative $(0,1)$, and such that any Type A or E vertex with an  incoming edge of derivative $(1,k)$ has an outgoing edge with derivative $(1,k+1)$. For such tropical curves $\gamma$, $\Aut_0\gamma$ is trivial, and $n_{\gamma_v}=1$ for every vertex apart from the Type C vertex $v_0$. Moreover, for such a curve, 
\[n_\gamma=n_{\gamma_{v_0}}=n_\Gamma\ .\]
Once $\Gamma=\{(1,1-m_1),\dotsc,(1,1-m_k)\}$ is chosen, such a tropical curve $\gamma$ with an isomorphism of $\gamma_{v_0}$ with $\Gamma$ is determined by choosing which $m_i$ Type E corners the $i$th edge exiting $v_0$ is connected to. Summing over  the possibilities  and using Lemma \ref{gamma homology} gives
\[\sum_{\gamma\vert \gamma_{v_0}=\Gamma}\frac{n_\gamma}{\abs{\Aut_0\gamma}}x^{g(\gamma)-1}q^{\beta(\gamma)}=n_\Gamma \prod_{i=1}^kx^{m_i-3}q^H\sigma_{m_i}=n_\Gamma\Psi(\Gamma)\]
where $\sigma_{m_i}$ is the $m_i$th elementary symmetric function in the $q^{-E_j}$. Plugging this into equation (\ref{GnFn}) gives the required formula (\ref{GnFn}) for $G_n$.

\stop

\

\section{Relative Gromov--Witten invariants of $\mathbb CP^{2}$ blown up at one point}\label{last}

In this section, we  calculate the relative Gromov--Witten invariants required in the previous section for Type E vertices. The required invariants are the relative Gromov--Witten invariants of $\mathbb CP^{2}$ blown up at $1$ point, relative to two lines $L_{1}$, $L_{2}$.

\begin{figure}[h]
\includegraphics{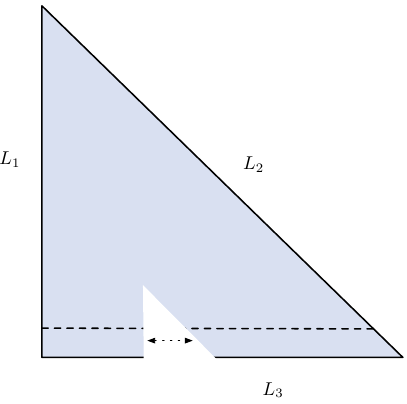}
\caption{}\label{tec14}
\end{figure}

To  picture $\mathbb CP^{2}$ blown up at a point, we can use Symington's almost toric blowup on $\mathbb CP^2$ to obtain an almost toric structure on our blowup; see section 4.3 of \cite{almosttoric}. Such an almost toric blowup produces a singular Lagrangian torus fibration with base depicted in Figure \ref{tec14}. The size of the little removed triangle represents the size of the symplectic ball removed to do a symplectic blowup; the remaining polytope should be regarded as glued along the two faces of this little removed triangle so that it has an integral-affine structure with a singularity at the top point of the little removed triangle. This singularity in the integral-affine structure reflects a focus-focus singularity in the Lagrangian torus fibration above it --- in other words, the torus fiber pinches to become a sphere which intersects itself once above this point. There are also elliptic singularities along the $3$ edges of this picture, as is usual for moment-map pictures. We are interested in $4$ holomorphic spheres in this picture. Over the left and righthand  boundaries are spheres $L_{1}$ and $L_{2}$ that are lines from $\mathbb CP^{2}$. Running down the glued-together edges of the little removed triangle is the exceptional sphere, and over the bottom boundary is a sphere $L_{3}$ that is the strict transform of a line passing through the point we blew up. 

We need the relative Gromov--Witten invariants of this space relative to $L_{1}$ and $L_{2}$, because Lemma \ref{Evertex} tells us that $\ex B\tc v$ is isomorphic to the explosion of this space relative to $L_1\cup L_2$ whenever $v$ is a Type E vertex.

We calculate these relative invariants  tropically by making a degeneration of this space into the two pieces above and below the dotted line. There is a connected family of exploded manifolds with one fiber isomorphic to $\ex B\tc v$, and another fiber an exploded manifold $\ex A'$ with smooth part the union of the two almost-toric manifolds above and below the dotted line in Figure \ref{tec14}, and stratification determined by $L_1$, $L_2$ and the dotted line in this figure. For pictures of tropical curves in $\totb{\ex A'}$, see figures \ref{tec21} and \ref{tec20}. The tropical completion of $\ex A'$ at the stratum with smooth part the bottom piece of Figure \ref{tec14} is  $\expl(\mathbb CP^{1},{0,\infty})\times\expl(\mathbb CP^1,0)$, which is a refinement of $\ex T\times\expl(\mathbb CP^1,0)$, so Lemma \ref{boundary count} suffices to determine its Gromov--Witten invariants.  The tropical completion of $\ex A'$ at the stratum with smooth part the top piece of Figure \ref{tec14} is the explosion of  the blowup of $\mathbb CP^{2}$ at a point,  relative to $L_{1}$, $L_{2}$, and $L_{3}$. Use $\ex A$ for this explosion.

\subsection{Gromov--Witten invariants of $\ex A$, or relative invariants of the blowup of $\mathbb CP^2$ at one point, relative to $L_1$, $L_2$, and $L_3$}

\

\begin{figure}[h]
\includegraphics{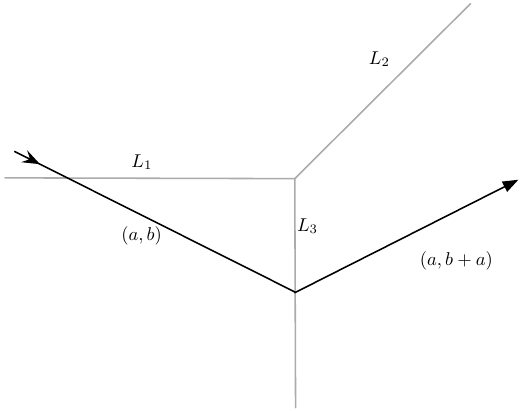}
\caption{A tropical curve in the tropical part of $\ex A$}\label{tec15}
\end{figure}

Figure \ref{tec15} shows the tropical part of the explosion, $\ex A$, of the one-point-blowup of $\mathbb CP^2$ relative to $L_{1}\cup L_{2}\cup L_{3}$, with a tropical curve.

To translate tropical curves in Figure \ref{tec15} to contact data with $L_{i}$,  unconstrained contact with $L_{1}$, $L_{2}$ or $L_{3}$ corresponds to an end with outgoing derivative $(-1,0)$, $(1,1)$, or $(0,-1)$ respectively. 

\begin{lemma}The contact data of curves in $\ex A$ obeys a balancing condition: the sum of the derivative of all ends (oriented outgoing) is some multiple of $(0,1)$. 
\end{lemma}
\pf This follows from the observation that $L_1$ and $L_2$ represent the same homology class in the manifold $\totl{\ex A}$. 

\stop

 The first Chern class of $\ex A$ vanishes, so the real virtual dimension of the  moduli space of curves in $\ex A$ with genus $g$ and unconstrained contact data $\Gamma$ is determined by $g$ and the number $k$ of ends of $\Gamma$.  
\[\text{virtual dimension}=2(k+g-1)\]
So, rigid curves either have genus $0$ and one end, or no ends and genus $1$. The homology class of a curve is determined by its contact data with $L_{i}$, and in particular, any curve in $\ex A$ without ends must be constant. Constant curves with genus $1$ are not stable, so we are left with rigid curves that are genus $0$ and have one end.  Our balancing condition then implies that all rigid holomorphic curves must be spheres with contact data consisting of one outgoing edge with derivative $(0,d)$. 

\

 To calculate the contribution $n_\gamma$ of a tropical curve $\gamma$ in $\totb{\ex A}$ to the Gromov--Witten invariants of $\ex A$, we need to know what relative Gromov--Witten invariants $n_{\gamma_v}$ to associate to each vertex. A vertex $v$ sent to the origin in $\totb{\ex A}$ uses the Gromov--Witten invariants of $\ex A$, because $\ex A\tc v=\ex A$. A vertex sent to  a two-dimensional stratum in $\totb{\ex A}$ uses the Gromov--Witten invariants of $\ex T^{2}$, because $\ex A\tc v=\ex T^2$. We discuss such Gromov--Witten invariants in Lemmas \ref{tropical dimension} and \ref{tropical count}. For $v$ on a one dimensional stratum $\ex A\tc v$ is a refinement of $\ex T^2$, however the tropical part of the refinement map $\ex A\tc v\longrightarrow \ex T^2$ is interesting for $v$ on the ray spanned by $(0,-1)$.
 
 \begin{lemma}\label{1ds} For $v$ on a $1$-dimensional stratum, $\ex A\tc v$ is a refinement of $\ex T^2$, isomorphic to $\ex T\times \expl (\ex CP^1,\{0,\infty\})$. The tropical part of the refinement map $\ex A\tc v\longrightarrow \ex T^2$ is as follows:
 \begin{itemize}
 \item for $v$ on the rays spanned by $(-1,0)$ or $(1,1)$, the tropical part of the refinement map is $v+x\mapsto x$;
 \item for $v$ on the ray spanned by $(0,-1)$, the tropical part of the refinement map is 
 \[v+(a,b)\mapsto (a,b-\max\{a,0\})\ .\]
 \end{itemize}
 \end{lemma}
 
 \proof

 A neighborhood, $N_i$ of $L_i$ in  $\totl{\ex A}$ is isomorphic to a neighborhood $U_i$ of a toric boundary component of a toric manifold. For $L_1$ and $L_2$, we can use neighborhoods of boundary strata of $\mathbb CP^2$, and for $L_3$, we can choose $U_3$ as a neighborhood of $\{0\}\times \mathbb CP^1$ in $(\mathbb CP^1)^2$. The image of the tropical part of $\expl N_i$ in Figure \ref{tec15} consists of the closure of the quadrants adjacent to the ray labeled $L_i$. There are obvious choices of isomorphism so that the tropical parts of the isomorphisms $\expl N_1\longrightarrow \expl U_1$ and $\expl N_2\longrightarrow \expl U_2$ are the identity, and the tropical part of the isomorphism $\expl N_3\longrightarrow \expl U_3$ is the map
  \[v+(a,b)\mapsto (a,b-\max\{a,0\})\ .\]
For $v$ in a $1$ dimensional stratum, $\ex A\tc v=\expl N_i\tc v=\expl U_i\tc v$, and $\expl U_i\tc v$ is the tropical completion of the explosion of the corresponding toric manifold at $v$. The explosion of any compact $2$-dimensional toric manifold is a refinement of $\ex T^2$, and the refinement map has tropical part the identity. The same holds for the tropical completion of such exploded manifolds at any point on a $1$-dimensional tropical stratum, and this tropical completion is isomorphic to $\ex T\times \expl(\mathbb CP^1,\{0,\infty\})$. Therefore, $\ex A\tc v$ is also a refinement of $\ex T^2$, isomorphic to $\ex T\times \expl(\mathbb CP^1,\{0,\infty\})$, and the tropical part of the refiment map $\ex A\tc v\longrightarrow \ex T^2$ is as specified above.

\stop

So, $\ex A\tc v$ is isomorphic to a refinement of $\ex T^2$ for any point $v$ in $\totb{\ex A}\setminus 0$. This implies that tropical curves in $\totb{\ex A}$ obey the usual balancing condition away from $0$ if we give $\totb{\ex A}$ the singular integral-affine structure obtained by cutting the plane along the line corresponding to $L_{3}$, and gluing the left and righthand sides of this cut so that a vector $(a,b)$ on the left side corresponds to the vector $(a,b+a)$ on the right.  In particular, the tropical curve of Figure \ref{tec15} corresponds to a `straight line'.
 
 We can now easily compute some Gromov--Witten invariants of $\ex A$. Let $n_{v_{1},\dotsc,v_{n};w}$ indicate the Gromov--Witten invariant of $\ex A$ that counts zero genus curves with contact data consisting of  incoming (constrained) ends with derivative $v_{i}$ and an outgoing (unconstrained) end with derivative $w$.
 
\begin{example}Consider computing  $n_{(0,1);(1,1)}$. The top tropical curve in Figure \ref{tec16} is the unique curve contributing to $n_{(0,1);(1,1)}$ when the incoming edge is constrained in the upper half plane. The bottom tropical curve is the only curve that contributes when the incoming edge is constrained in the lower half plane.
 
\begin{figure}[h] 
\includegraphics{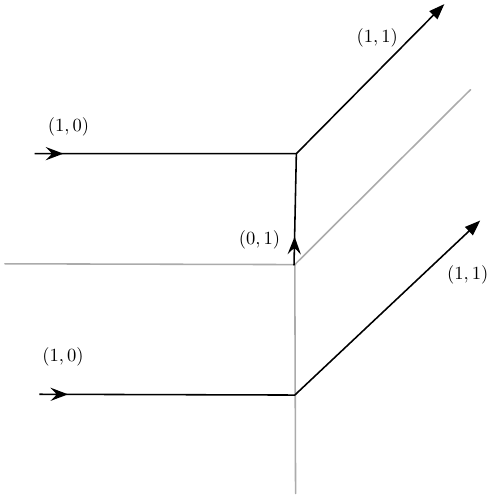}
\caption{}\label{tec16}
\end{figure}

\noindent Therefore, using our simplified gluing formula from Theorem \ref{simplified gluing},
\[n_{(1,0);(1,1)}=1=n_{;(0,1)}\ .\]

\end{example}

\begin{lemma} $n_{(a,b);(a,b+d)}$ depends only on $a$ and $d$, so
\[n_{(a,b+1);(a,c)}=n_{(a,b);(a,c-1)}\ .\]

\end{lemma}

\pf 

\begin{figure}[h]
\includegraphics{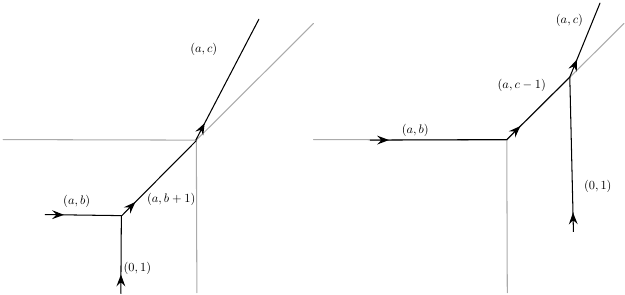}
\caption{}\label{tec17}
\end{figure}

For $a>0$, we can equate the two ways of calculating $n_{(a,b),(0,1);(a,c)}$ shown in Figure \ref{tec17}. 

\noindent So
\[n_{(a,b+1);(a,c)}=n_{(a,b);(a,c-1)}\]
as required.  

\stop

Calculating $n_{(a,b);(a,b+a)}$ by constraining the incoming edge below the ray $(-a,-b)$ gives that 
\begin{equation}n_{(a,b);(a,b+a)}=1\end{equation}
and
\begin{equation}n_{(a,b);(a,b+d)}=0\text{ if }d>\abs a\ .\end{equation}

Similarly, calculating constraining the incoming edge above the ray $(-a,-b)$ gives that
\begin{equation}n_{(a,b);(a,b+d)}=0 \text{ if }d<0\ .\end{equation}

\begin{prop} \label{abmult}\[n_{(a,b);(a,b+d)}={a\choose d}\]

\end{prop}

\begin{figure}[h]
\includegraphics{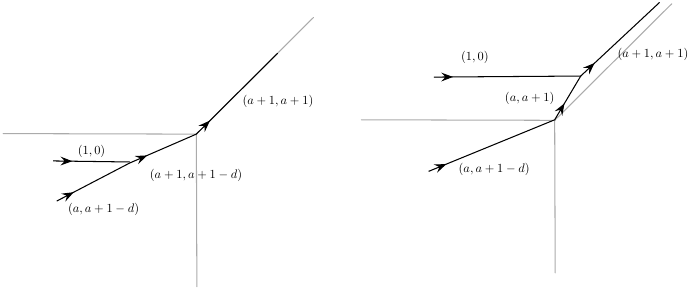}
\caption{}\label{tec18}\end{figure}

\pf
If $0\leq d\leq a$, we can equate the two ways of calculating $n_{(a,a+1-d),(1,0);(a+1,a+1)}$ shown in Figure \ref{tec18}.

\noindent Therefore 
\[(a+1-d)n_{(a+1,a+1-d);(a+1,a+1)}=(a+1)n_{(a,a+1-d);(a,a+1)}\ .\]
Starting with the case 
\[n_{(a,b);(a,b+a)}=1={a\choose a}\]
induction on $d\geq a$ using the above equation gives 
\[n_{(a,b);(a,b+d)}={a\choose d} \] 
as required. 

\stop

\

We still need to compute $n_{;(0,d)}$.

\begin{lemma}\label{dmult} \[n_{;(0,d)}=\frac{(-1)^{d+1}}{d^2}\]
\end{lemma}
 
\begin{figure}[h]
\includegraphics{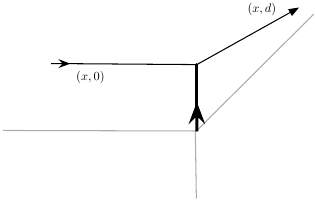}
\caption{}\label{tec19}\end{figure}

 \pf
 Consider $p(x)=n_{(x,0);(x,d)}$. For $x\geq 0$, we can compute $p(x)$ by restricting the incoming edge $(x,0)$ in the upper half plane. The gluing formula for $p(x)$ uses tropical curves in the form shown in Figure \ref{tec19}. In this diagram, the thick edge indicates some number of edges with upward pointing derivatives adding to $(0,d)$.

\noindent A diagram with $k$ such upward pointing edges contributes some fixed multiple of $x^{k}$ to $p(x)$, so $p(x)$ is a polynomial in $x$ of degree at most $d$. The only possibility for a single such edge is when its derivative is $(0,d)$. This tropical curve contributes $xdn_{;(0,d)}$ to $p(x)$, therefore the coefficient of the linear term in $p(x)$ is $dn_{;(0,d)}$. We already know that $p(x)=0$ for $0\leq x< d$, and $p(d)=1$, therefore
\[p(x)=\frac 1{d!}\prod_{i=0}^{d-1}(x-i)\ .\]
The coefficient of the linear term of $p(x)/d$ then gives us $n_{;(0,d)}$.
\[n_{;(0,d)}=\frac1{d!d}\prod_{i=1}^{d-1}(-i)=\frac{(-1)^{d-1}}{d^{2}}\]
Although we have only argued that the above equation holds for $d\geq 1$, a very similar argument also gives that $n_{;(0,d)}=(-1)^{d+1}/d^{2}$ for $d<0$.
Of course, the equations determining $n_{;(0,d)}$ are massively overdetermined. A combinatorially-talented thinker could deduce the above equations for $d>0$ and the corresponding formula for $n_{(a,b);(a,b+d)}$  simply from knowing that $n_{;(0,1)}=1$ and $n_{(1,0);(1,d)}=0$ for $d>1$.
\stop

Because the curves contributing to $n_{;(0,d)}$ do not contact the divisor $L_3$, and the same is true for the corresponding curves in $\ex B\tc v$ for Type E vertices, Lemma \ref{dmult} implies Equation (\ref{n0d}). We also analyze the Gromov--Witten invariants of $\ex B\tc v$ using the Gromov--Witten invariants of  $\ex A'$ in the next section.

\subsection{Gromov--Witten invariants of $\ex A'$}

\begin{figure}[h]
\includegraphics{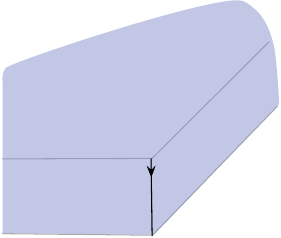}
\caption{The tropical part of the exceptional sphere in the tropical part of $\ex A'$}\label{tec21}
\end{figure}

The exploded manifold $\ex A'$ is in a connected family of exploded manifolds also containing the blowup of $\mathbb CP^2$ relative to $L_1$ and $L_2$. We need to analyze tropical curves in the tropical part of $\ex A'$, pictured in Figure \ref{tec21}. 

For calculating the contribution of a tropical curve $\gamma$ in $\totb{\ex A'}$ to Gromov--Witten invariants of $\ex  A'$, vertices above the bottom boundaries of $\totb{\ex A'}$ in Figure \ref{tec21} contribute the same as the corresponding vertices in $\totb{\ex A}$. In particular, 
\begin{itemize}
\item for $v$ at the interior $0$--dimensional stratum, $\ex A'\tc v=\ex A$;
\item for $v$ in a $2$--dimensional stratum, $\ex A'\tc v=\ex T^2$;
\item for $v$ in any interior $1$--dimensional stratum, $\ex A'\tc v$ is a refinment of $\ex T^2$, and the tropical part of the refinment map is as described in Lemma \ref{1ds}.
\end{itemize}

We also have
\begin{lemma}\label{lb}
For $v$ at the $0$--dimensional stratum on the boundary of $\totb{\ex A'}$, the tropical completion $\ex A'\tc v$  is isomorphic to $\expl(\mathbb CP^1,\{0,\infty\})\times\expl(\mathbb CP^1,0)$, and the tropical part of the isomorphism is
\[v+(a,b)\mapsto (a,b-\max{a,0}) \ .\]

If $v$ is on either of the boundary $1$--dimensional strata,  $\ex A'\tc v$ is isomorphic to $\ex T\times \expl(\mathbb CP^1,0)$, and the tropical part of this isomorphism sends $v+(0,b)$ to $(0,b)$.
\end{lemma}

\pf  For  $v$ at the boundary $0$--dimensional stratum of $\totb{\ex A'}$, $\ex A'\tc v$ is the explosion of the closure of the smooth part of that stratum, depicted in Figure \ref{tec14} as the region below the dotted line. This manifold is isomorphic to $(\mathbb CP^1)^2$, with the divisor from the dotted line and the left and right hand boundaries, so $\totl{\ex A'\tc v}$ is isomorphic to $(\mathbb CP^1,\{0,\infty\})\times(\mathbb CP^1,0)$, and $\ex A'\tc v=\expl(\mathbb CP^1,\{0,\infty\})\times\expl(\mathbb CP^1,0)$, as required. 

If $v'$ is a stratum of $\totb{\ex A'}$ with closure containing $v$, there is a natural identification of  $\ex A'\tc {v'}$ with $(\ex A'\tc v)\tc{v'}$. It follows that for $v'$ in any boundary $1$--dimensional stratum $\ex A'\tc {v'}$ is the tropical completion of $\expl(\mathbb CP^1,\{0,\infty\})\times \expl(\mathbb CP^1,0)$ at a stratum corresponding to $0$ or $\infty$ in the first $\mathbb CP^1$. In particular, $\ex A'\tc{v'}$ is isomorphic to $\ex T\times\expl(\mathbb CP^1,0)$. The tropical part of these isomorphisms must preserve the direction corresponding to the divisor $0$ in the second $\mathbb CP^1$, and this must be $(0,1)$ because it coincides with the direction of the vertical $1$--dimensional stratum.  

Similarly,  using the natural identifications of $\ex A'\tc {v'}$ with $(\ex A'\tc v)\tc {v'}$ when the closure of the $v'$--stratum of $\totb{\ex A'}$ contains $v$, we can show that the $1$--dimensional strata that appear parallel in Figure \ref{tec21} are actually parallel. It follows that the tropical part of our isomorphisms $\ex A'\tc v\longrightarrow \expl(\mathbb CP^1,\{0,\infty\})\times \expl(\mathbb CP^1,0)$ is 
 \[v+(a,b)\mapsto (a,b-\max{a,0}) \]
 as required.

\stop  

With the above, Lemma \ref{boundary count}, and the observation that $\expl(\mathbb CP^1,\{0,\infty\})\times \expl(\mathbb CP^1,0)$ is a refinement of $\ex T\times \expl(\mathbb CP^1,0)$, we can now determine the Gromov--Witten invariants of $\ex A'$ from the Gromov--Witten invariants of $\ex A$. 

\begin{lemma}\label{AA'} Suppose that $\gamma$ is a rigidly constrained tropical curve in $\totb{\ex A'}$ with $n_\gamma\neq 0$, and any constrained ends constrained  within the interior of $\totb{\ex A'}$. Then, $\gamma$ touches the boundary of $\totb{\ex A}$ at at most one point, and at this point, $\gamma$ has a univalent vertex with incoming derivative $(0,-1)$. Let $\gamma'$ be the tropical curve in $\totb{\ex A}$ obtained from $\gamma$ by removing the univalent vertex, and extending the edge to be infinite. Then
\[n_\gamma=n_\gamma'\ .\]
where $n_\gamma'$ counts curves in $\ex A$ with the same ends constrained, (so the extra end is not constrained.) 
\end{lemma}

\pf

A dimension count, our tropical balancing conditions and Lemma \ref{dimension vanishing} imply  that we have $3$ cases:
\begin{description}
\item[Case 1]  All but one end of $\gamma$ is constrained, $n_\gamma$ counts zero genus curves, and $\gamma$ does not touch the boundary of $\totb{\ex A'}$.
\item[Case 2] All ends of $\gamma$ are constrained, $n_\gamma$ counts genus one curves, and $\gamma$ does not touch the boundary of $\totb{\ex A'}$
\item[Case 3] All ends of $\gamma$ are constrained, and $\gamma$ has a unique vertex on the boundary of $\totb{\ex A'}$, and this vertex has incoming derivative $(0,-1)$. In this case $n_\gamma $ counts zero genus curves, and in particular $\gamma$ has zero genus.
\end{description}

In cases $1$ and $2$, it is clear that $n_\gamma=n_{\gamma'}$. It remains to prove this equality in Case 3.

In this case, because $\gamma$ has zero genus, we can apply Theorem \ref{simplified gluing} with edges of $\gamma$ oriented so that each vertex on the interior has a unique outgoing edge, and the boundary vertex has its edge oriented inwards. For this boundary vertex $v$, $n_{\gamma_v}=1$, by Lemmas \ref{lb} and \ref{boundary count}. At other vertices $v'$, $n_{\gamma'_{v'}}=n_{\gamma_{v'}}$, so Theorem \ref{simplified gluing} implies $n_\gamma=n_{\gamma'}$, as required.

\stop

\begin{figure}[h]
\includegraphics{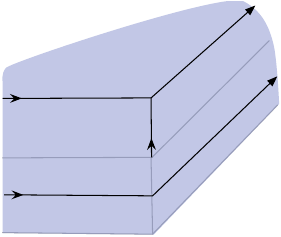}
\caption{}\label{tec20}\end{figure}

 For example, the contribution of the tropical curve in Figure \ref{tec21} is $1$ (corresponding to the exceptional curve). Either of the tropical curves in Figure \ref{tec20} may be used to calculate that $1$ is the Gromov--Witten invariant of $\ex A'$ with contact data a constrained edge entering with derivative $(1,0)$ and one unconstrained edge exiting with  derivative  $(1,1)$. Moreover, this counts genus $0$ curves that intersect the exceptional sphere once.
 
\begin{figure}[h]
\includegraphics{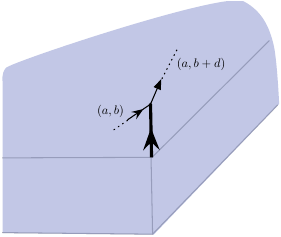}
\caption{}\label{tec22}
\end{figure}

As $n_{;(0,d)}=(-1)^{d+1}/d^{2}$ implies Proposition \ref{abmult},  whenever we see a part of a tropical curve $\gamma$ looking like Figure \ref{tec22} --- where the thick edge may be replaced by many edges with derivatives adding up to $(0,d)$ --- the total effect of this part of $\gamma$ is to multiply by $a\choose d$, and to affect the corresponding homology class  by adding $d$ to its intersection with the exceptional divisor. In the tropical part,  $\totb{\ex B}$, of the exploded manifold we used to represent $\mathbb CP^{2}$ blown up at $n$ points,  we may intuitively understand this as saying that a rigid edge with derivative $(a,b)$ may interact with rigid edges coming up from the $i$th lower corner to leave with derivative $(a,b+d)$. After summing over all possibilities, this introduces a factor $a\choose d$ to the Gromov--Witten invariant, and corresponds to intersecting the $i$th exceptional sphere $d$ times.

\bibliographystyle{plain}
\bibliography{ref.bib}

\begin{thebibliography}{10}

\bibitem{acgw}
Dan Abramovich and Qile Chen.
\newblock Stable logarithmic maps to {D}eligne-{F}altings pairs {II}.
\newblock {\em The Asian Journal of Mathematics}, 18(3):465--488, 2014.

\bibitem{acgsdegeneration}
Dan Abramovich, Qile Chen, Mark Gross, and Bernd Siebert.
\newblock Decomposition of degenerate {G}romov--{W}itten invariants.
\newblock \href{https://arxiv.org/abs/1709.09864}{arXiv:1709.09864}.

\bibitem{ilgw}
Dan Abramovich and Jonathan Wise.
\newblock Birational invariance in logarithmic {G}romov-{W}itten theory.
\newblock {\em Compos. Math.}, 154, 2018.

\bibitem{Brugalle}
Erwan Brugall\'e.
\newblock Floor diagrams relative to a conic, and {GW}-{W} invariants of del
  {P}ezzo surfaces.
\newblock {\em Adv. Math.}, 279:438--500, 2015.

\bibitem{BryanLeung}
Jim Bryan and Naichung~Conan Leung.
\newblock The enumerative geometry of {$K3$} surfaces and modular forms.
\newblock {\em J. Amer. Math. Soc.}, 13(2):371--410 (electronic), 2000.

\bibitem{CH}
Lucia Caporaso and Joe Harris.
\newblock Counting plane curves of any genus.
\newblock {\em Invent. Math.}, 131(2):345--392, 1998.

\bibitem{almosttoric}
Naichung Conan~Leung and Margaret Symington.
\newblock Almost toric symplectic four-manifolds.
\newblock {\em Journal of Symplectic Geometry}, 8, 01 2004.

\bibitem{GPqc}
L.~G{\"o}ttsche and R.~Pandharipande.
\newblock The quantum cohomology of blow-ups of {${\bf P}^2$} and enumerative
  geometry.
\newblock {\em J. Differential Geom.}, 48(1):61--90, 1998.

\bibitem{GSlogGW}
Mark Gross and Bernd Siebert.
\newblock Logarithmic {G}romov-{W}itten invariants.
\newblock {\em J. Amer. Math. Soc.}, 26(2):451--510, 2013.

\bibitem{IonelGW}
Eleny-Nicoleta Ionel.
\newblock G{W} invariants relative to normal crossing divisors.
\newblock {\em Adv. Math.}, 281:40--141, 2015.

\bibitem{IP}
Eleny-Nicoleta Ionel and Thomas~H. Parker.
\newblock The symplectic sum formula for {G}romov-{W}itten invariants.
\newblock {\em Ann. of Math. (2)}, 159(3):935--1025, 2004.

\bibitem{scgp}
Brett Parker.
\newblock Notes on exploded manifolds and a tropical gluing formula for
  {G}romov-{W}itten invariants.
\newblock \href{http://arxiv.org/abs/1605.00577}{arXiv:1605.00577}.

\bibitem{egw}
Brett Parker.
\newblock Gromov-{W}itten invariants of exploded manifolds.
\newblock \href{http://arxiv.org/abs/1102.0158}{arXiv:1102.0158}, 2011.

\bibitem{iec}
Brett Parker.
\newblock Exploded manifolds.
\newblock {\em Adv. Math.}, 229:3256--3319, 2012.
\newblock \href{http://arxiv.org/abs/0910.4201}{arXiv:0910.4201}.

\bibitem{elc}
Brett Parker.
\newblock Log geometry and exploded manifolds.
\newblock {\em Abh. Math. Sem. Hamburg}, 82:43--81, 2012.
\newblock \href{http://arxiv.org/abs/1108.3713}{arxiv:1108.3713}.

\bibitem{evc}
Brett Parker.
\newblock Holomorphic curves in exploded manifolds: Kuranishi structure.
\newblock \href{http://arxiv.org/abs/1301.4748}{arXiv:1301.4748}, 2013.

\bibitem{tropicalIonel}
Brett Parker.
\newblock On the value of thinking tropically to understand {I}onel's {GW}
  invariants relative normal crossing divisors.
\newblock \href{http://arxiv.org/abs/1407.3020}{arXiv:1407.3020}, 2014.

\bibitem{cem}
Brett Parker.
\newblock Holomorphic curves in exploded manifolds: compactness.
\newblock {\em Adv. Math.}, 283:377--457, 2015.
\newblock \href{http://arxiv.org/abs/0911.2241}{arXiv:0911.2241}.

\bibitem{gfgw}
Brett Parker.
\newblock Tropical gluing formulae for {G}romov-{W}itten invariants.
\newblock \href{http://arxiv.org/abs/1703.05433}{arXiv:1703.05433}, 2017.

\bibitem{dre}
Brett Parker.
\newblock De {R}ham theory of exploded manifolds.
\newblock {\em Geometry and Topology}, 22(1):1--54, 2018.
\newblock \href{http://arxiv.org/abs/1003.1977}{arXiv:1003.1977}.

\bibitem{vfc}
Brett Parker.
\newblock Holomorphic curves in exploded manifolds: virtual fundamental class.
\newblock {\em Geometry and Topology}, 23:1877--1960, 2019.
\newblock \href{http://arxiv.org/abs/1512.05823}{arXiv:1512.05823}.

\bibitem{rang}
Dhruv Ranganathan.
\newblock Logarithmic {G}romov-{W}itten theory with expansions.
\newblock {\em Algebr. Geom.}, 9(6):714--761, 2022.

\bibitem{Glog}
Dhruv Ranganathan and Jonathan Wise.
\newblock Rational curves in the logarithmic multiplicative group.
\newblock {\em Proc. Am. Math. Soc.}, 148(1):103--110, 2020.

\bibitem{SS}
Mendy Shoval and Eugenii Shustin.
\newblock On {G}romov-{W}itten invariants of del {P}ezzo surfaces.
\newblock {\em Internat. J. Math.}, 24(7):1350054, 44, 2013.

\bibitem{VakilSurfaces}
Ravi Vakil.
\newblock Counting curves on rational surfaces.
\newblock {\em Manuscripta Math.}, 102(1):53--84, 2000.

\end{thebibliography}

\end{document}